\documentclass[12pt]{amsart}

\usepackage{amsfonts,amsmath,amssymb}
\usepackage{latexsym,graphicx}
\usepackage{chicago}
\usepackage{graphicx,multirow,rotating}

\newtheorem{defi}{Definition}[section]

\newtheorem{theo}{Theorem}[section]
\newtheorem{lem}{ Lemma}[section]

\newtheorem{remark}{Remark}[section]

\newcommand{\eref}[1]{(\ref{#1})}
\newcommand\ind{{{{1}}\hspace{-1,1mm}{\mathrm I}}}

\newcommand{\pr}[0]{\mathrm{I}\!\mathrm{P}}

\newcounter{hypc}

\newcommand{\hypothese}[1]{\stepcounter{hypc}\tag{$\mathbf{A_{\thehypc}}$}
\label{#1}}

\newcounter{regc}

\newcounter{noisc}
\newcommand{\noise}[1]{\stepcounter{noisc}\tag{$\mathbf{N_{\thenoisc}}$}
\label{#1}}

\newcounter{identc}
\newcommand{\ident}[1]{\stepcounter{identc}\tag{$\mathbf{I1_{\theidentc}}$}
\label{#1}}

\newcounter{identtc}

\newcounter{condc}
\newcommand{\condition}[1]{\stepcounter{condc}\tag{$\mathbf{C_{\thecondc}}$}
\label{#1}}

\def\build#1_#2^#3{\mathrel{\mathop{\kern 0pt#1}\limits_{#2}^{#3}}}
\def\cvl{\build{\ \longrightarrow\ }_{\nti}^{{\mathcal L}}}

\def\cvp{\build{\ \longrightarrow\ }_{\nti}^{{\pr}}}
\def\cvL1{\build{\ \longrightarrow\ }_{\nti}^{{\mathbb{L}^1}}}
\def\cv{\build{\ \longrightarrow\ }_{\nti}^{}}
\def\nti{n \to \infty}

\title{Estimation in autoregressive model with measurement error}
\author{   J\'er\^ome    DEDECKER$^{(1)}$,   Adeline   Samson$^{(1)}$,
  Marie-Luce TAUPIN$^{(2)}$}

\begin{document}

\maketitle

\noindent{\small
(1)  Laboratoire MAP5  UMR  CNRS 8145,  Universit\'e Paris  Descartes,
Sorbonne Paris Cit\'e,\\
(2) Laboratoire Statistique et G\'enome UMR CNRS 8071- USC INRA,
Universit\'e d'\'Evry Val d'Essonne}


\begin{abstract}
Consider an autoregressive model
with measurement error: we observe $Z_i=X_i+\varepsilon_i$, where
$X_i$ is a stationary solution of the autoregressive equation  $X_i=f_{\theta^0}(X_{i-1})+\xi_i$. The regression function
$f_{\theta^0}$ is known up to a finite dimensional parameter $\theta^0$. The distributions of
 $X_0$ and  $\xi_1$ are unknown
whereas the distribution  of $\varepsilon_0$ is completely known.
We want to estimate the parameter $\theta^0$ by using the observations
$Z_0,\ldots,Z_n$. We propose an estimation procedure
based on a modified least square criterion.
This procedure provides an asymptotically normal estimator $\hat \theta$ of $\theta^0$, for a large
class of regression functions and  various noise distributions.
\end{abstract}

\medskip

{\small

\noindent {\bf AMS 2000 MSC}: Primary 62J02, 62F12, Secondary 62G05, 62G20.

\section{Introduction}
\setcounter{equation}{0}
\setcounter{lem}{0}
\setcounter{theo}{0}
We consider an autoregressive model with measurement error satisfying
\begin{eqnarray}
\label{model}
\left\lbrace
\begin{array}{ll}
Z_i&=X_i+\varepsilon_i, \\
X_i&=f_{\theta^0}(X_{i-1})+\xi_i
\end{array}
\right.
\end{eqnarray}
where one
observes $Z_0,\cdots,Z_n$ and the random variables $\xi_i,X_i,\varepsilon_i$ are unobserved.
The regression function $f_{\theta^0}$ is known up to a finite dimensional
parameter
$\theta^0$, belonging to the interior $\Theta^\circ$ of a
compact set  $\Theta \subset\mathbb{ R}^{d}$.
The centered innovations $(\xi_i)_{i\geq 1}$ and the errors
$(\varepsilon_i)_{i\geq 0}$ are  independent and
identically   distributed  (i.i.d.)   random  variables   with  finite
variances
$\mbox{Var}(\xi_1)=\sigma_{\xi}^2$ and  $\mbox{Var}(\varepsilon_0)=
\sigma_{\varepsilon}^2$.
We assume that $\varepsilon_0$  admits
 a known density with respect to the Lebesgue measure, denoted by
$f_{\varepsilon}$. Furthermore we assume that the random variables $X_0$, $(\xi_i)_{i \geq 1}$ and  $(\varepsilon_i)_{i \geq 0}$ are
independent.
The distribution
of  $\xi_1$ is unknown  and does  not necessarily  admit a  density with
respect  to the  Lebesgue measure.   We  assume that  $(X_i)_{i \geq 0}$ is  strictly
stationary,
which means that the initial distribution of $X_0$ is an invariant distribution for the transition kernel
of the homogeneous Markov chain $(X_i)_{i \geq 0}$.

Our aim is to estimate $\theta^0$ for a large class of functions $f_\theta$, whatever  the known error distribution, and without the
knowledge of the $\xi_i$'s distribution.
The distribution of the  innovations being unknown, this model belongs
to the family of
semi-parametric models.

\subsection*{Previously known results}

Several authors have considered the case where the function $f_\theta$ is linear (in both $\theta$ and $x$), see \textit{e.g.}
Andersen and Deistler \citeyear{ANDDEIS}, Nowak \citeyear{Nowak},
Chanda \citeyear{Chanda,Chanda96},
  Staudenmayer  and
Buonaccorsi \citeyear{StaudenmayerBuonaccorsi}, and Costa   \textit{et  al.}   \citeyear{CostaAlpuim}.
We can note that, in this specific case, the model \eref{model} is also an ARMA model
(see Section \ref{Linear}  for further details). Consequently, all previously
known   estimation  procedures   for   ARMA  models   can  be   applied
here, without assuming that the error distribution is known.

For a general regression
function,  the model \eref{model} is  a Hidden Markov Model with
possibly a non
compact continuous
state space, and with unknown innovation distribution. When the innovation distribution is known up to a finite
dimensional parameter, the model \eref{model} is fully parametric and
various results are already stated.  Among others, the parameters can be estimated by
maximum   likelihood,  and   consistency,   asymptotic  normality   and
efficiency have been proved.
For further references on estimation in fully parametric Hidden Markov
Models, we refer for instance to
Leroux \citeyear{Leroux},
Bickel \textit{et al.} \citeyear{BickelRitovRyden}, Jensen and Petersen \citeyear{JensenPetersen},
Douc and Matias \citeyear{DoucMatias}, Douc \textit{ et al.} \citeyear{DoucMoulinesRyden}, Fuh \citeyear{Fuh2006},
Genon-Catalot  and   Laredo  \citeyear{VGCCL},  Na   \textit{et.  al.}
\citeyear{Naetal}, and Douc \textit{et al.} \citeyear{DoucMoulinesOlssonvanHandel}.

In this paper, we consider the case where the innovation distribution is unknown, and thus the model is
not fully parametric. In this general context,  there are few results. To our knowledge,
the only paper which gives a consistent estimator is the paper by
Comte and Taupin \citeyear{ComteTaupin2001}.  These authors propose an
estimation procedure based on a modified least squares minimization.
They  give an upper bound for the rate
of convergence of their estimator, that depends  on  the smoothness  of the
regression function  and on the  smoothness of $f_\varepsilon$.
  Those results are obtained by assuming that the distribution $P_X$ of $X_0$ admits
a density $f_X$  with respect  to the
Lebesgue measure and that
the stationary Markov
chain $(X_i)_{i\geq  0}$ is absolutely regular
($\beta$-mixing). The main drawback of their approach is that their estimation criterion is not explicit, hence the links
between  the convergence rate of their estimator and the smoothness  of the regression  function and of the  error distribution  are not explicit either.
Consequently, Comte and Taupin \citeyear{ComteTaupin2001} are able to prove that their estimator achieves the
parametric  rate only for  very few  couples of  regression functions/error
distribution.    Lastly   their   dependency  conditions   are   quite
restrictive, and the assumption that $X$ admits a density is not natural in this context.

\subsection*{Our results}

In this paper, we propose a new estimation procedure which provides a consistent
estimator with a parametric rate of convergence in a very general context.
Our approach is
based on the new contrast function
\begin{eqnarray*}
 S_{\theta^0,P_X}(\theta)=\mathbb{E}[(Z_1-f_{\theta}(X_{0}))^2 \, w(X_{0})],
\end{eqnarray*}
where $w$ is a weight function to be chosen and $\mathbb{E}$ is the expectation
$\mathbb{E}_{\theta^0,P_X}$. We assume that $w$ is such  that $(wf_\theta)^*/f_\varepsilon^*$ and
$(wf_\theta^2)^*/f_\varepsilon^*$ are integrable, where  $\varphi^*$ is the Fourier transform of a function $\varphi$.
We  estimate $\theta^0$ by
$
\widehat{\theta}=\arg\min_{\theta \in
\Theta}
S_{n}(\theta),\label{thetacintro}
$
where
\begin{equation}\label{Snintro}
 S_{n}(\theta)
=
\frac{1}{2\pi n}\sum_{k=1}^n
\mathbb{R}e
\int \frac{\left(\big(Z_k-f_{\theta}\big)^2w\right)^*(t)\; e^{-itZ_{k-1}}}{f_\varepsilon^*(-t)}dt,
\end{equation}
where $\mathbb{R}e(u)$ is the real part of $u$.
Under general assumptions, we prove that the estimator defined $\widehat{\theta}$ is consistent.
Moreover, we give some conditions under which the parametric rate
of   convergence  as  well   as  the   asymptotic  normality   can  be
stated.  Those results hold under weak dependency conditions
as
introduced in Dedecker and Prieur \citeyear{JDCPCoeff}.

This procedure is clearly simpler than that of Comte and Taupin
\citeyear{ComteTaupin2001}.  The resulting  rate is  more
explicit and  links directly the smoothness of  the regression function
to  that  of  $f_\varepsilon$.  Our new  estimator  is  asymptotically
Gaussian for a large class of regression
functions, which is not the case in
 Comte and Taupin \citeyear{ComteTaupin2001}.
%

The asymptotic properties of our estimator are illustrated through a
simulation  study. It  confirms  that our  estimator  performs well  in
various contexts, even in cases where the Markov chain $(X_i)_{i \geq 0}$
is not $\beta$-mixing (and not even irreducible), when the ratio signal to
noise is small  or large, for various sample  sizes, and for different
types of error distribution. Our estimator always better performs
than  the  so-called naive  estimator  (built  by replacing  the
non-observed $X$  by $Z$  in the usual  least squares  criterion). Our
estimation procedure depends on the choice of the weight function $w$. The influence
of this weight function is also studied in the
simulations.

 Finally, we propose a more general estimator when it is not possible to find a weight function $w$  such  that $(wf_\theta)^*/f_\varepsilon^*$ and
$(wf_\theta^2)^*/f_\varepsilon^*$ are integrable.
We establish a consistency
result,
and we give an upper bound for the quadratic risk, that relates the
smoothness  properties of  the  regression function  to that   of
$f_\varepsilon$.
These last results are proved under $\alpha$-mixing conditions.

The paper is organized as follows.
In Section \ref{estim} we present our estimation procedure. The
theoretical properties of the estimator are stated in Section \ref{proprietes}. The simulations are presented in Section \ref{simu_context}.
In Section \ref{generalcrit} we introduce a more  general estimator and we describe its asymptotic behavior. The proofs are gathered in Appendix.

\section{Estimation procedure}
\label{estim}
In  order to  define more  rigorously the  criterion presented  in the
introduction, we first give some preliminary notations and assumptions.

\subsection{Notations}
Let
$$\parallel \varphi\parallel_1=\int \vert \varphi(x)\vert
dx,\; \parallel \varphi\parallel_2^2=\int \varphi^2(x)dx, \;  \mbox{ and }\parallel
\varphi\parallel_\infty=\sup_{x\in \mathbb{R}}\vert \varphi(x)\vert.$$
The convolution product of two square integrable
functions $p$ and $q$ is denoted by
$p \star q(z) = \int p(z-x) q(x) dx$. The Fourier transform $\varphi^*$ of a function $\varphi$  is defined by
$$\varphi^*(t)=\int e^{itx}\varphi(x)dx.$$
For
$\theta \in \mathbb{R}^d$, let $\parallel \theta\parallel_{\ell^2}^2=\sum_{k=1}^d
\theta_k^2$, and let $\theta^\top$ be the transpose matrix of $\theta$.

For a map $(\theta,u) \mapsto \varphi_\theta(u)$ from $\Theta \times \mathbb{R}$ to $\mathbb{R}$,
 the first and
second derivatives with respect to $\theta$ are denoted by
\begin{eqnarray*}
\varphi^{(1)}_{\theta}(\cdot)&=&\left(\varphi_{\theta,j}^{(1)}(\cdot)\right)_{1\leq
  j\leq d}, \
\mbox{
with }\varphi_{\theta,j}^{(1)}(\cdot)=\frac{\partial \varphi_{\theta}(\cdot)}{\partial
  \theta_j}\mbox{ for
}j\in \{1,\cdots,d\}\\
\mbox{ and
}\hspace{0.5cm}\varphi^{(2)}_{\theta}(\cdot)&=&\left(\varphi_{\theta,j,k}^{(2)}(\cdot)\right)_{1\leq
  j,k\leq d}, \
\mbox{ with }\varphi_{\theta,j,k}^{(2)}(\cdot)=\frac{\partial^2
  \varphi_{\theta}(\cdot)}{\partial \theta_j \partial \theta_k},\mbox{ for
}j,k\in \{1,\cdots,d\}.
\end{eqnarray*}
From now, $\mathbb{P}$, $\mathbb{E}$ and $\mbox{Var}$ denote respectively the
probability $\mathbb{P}_{\theta^0,P_X}$, the expected value
$\mathbb{E}_{\theta^0,P_X}$ and  the variance
$\mbox{Var}_{\theta^0,P_X}$, when the underlying and unknown true
parameters are $\theta^0$ and $P_X$.

\subsection{Assumptions}

We consider three types of assumptions.

\bigskip

\noindent $\bullet$ \textbf{Smoothness and moment assumptions}
\begin{align}
\hypothese{df} &\mbox{On $\Theta^\circ$, the function
}\theta\mapsto f_\theta
\mbox{ admits continuous derivatives with respect to } \theta \mbox{ up to the }\\
  &\notag  \mbox{order 3}.
\\
&     \hypothese{emotilde}    \mbox{On     $\Theta^\circ$,    the
  quantity } w(X_{0})(Z_1-f_{\theta}(X_{0}))^{2},
\mbox{ and the absolute values of its derivatives} \\
&\notag\mbox{with respect to } \theta
 \mbox{ up to order 2  have a finite expectation.}
\end{align}

\noindent $\bullet$ \textbf{Identifiability assumptions}
\begin{align}
&\ident{Spttilde}\mbox{ The quantity }
S_{\theta^0,P_X}(\theta)=
\mathbb{E}[(f_{\theta^0}(X)-f_{\theta}(X))^2w(X)] \mbox{ admits one unique minimum at}\\&\notag
\theta=\theta^{0}.
\\
& \ident{emo2tilde}\mbox{ For all }\theta\in\Theta^\circ,  \mbox{ the matrix }
 S_{\theta^0,P_X}^{(2)}(\theta)=\left(
\frac{\partial^2   S_{\theta^0,P_X}(\theta)}{\partial                  \theta_i\partial
   \theta_j}\right)_{1\leq i,j\leq d} \mbox{ exists and the matrix }\\&\notag
  S^{(2)}_{\theta^0,P_X}(\theta^0)
=2\, \mathbb{E}\left[w(X)\left( f^{(1)}_{\theta^0}(X)\right)\left( f^{(1)}_{\theta^0}(X)
\right)^\top\right]\mbox{ is positive definite}.
\end{align}

\noindent $\bullet$ \textbf{Assumptions on $f_\varepsilon$}
\begin{align}
\noise{fepsnn}&\mbox{The density } f_{\varepsilon} \mbox{ belongs to
  }\mathbb{L}_2(\mathbb{R})  \mbox{ and for all }
x\in \mathbb{R},\;f_{\varepsilon}^*(x)\not=0.
\end{align}

\medskip

The   assumption  \eref{fepsnn}  is   quite  usual   when  considering
estimation in the convolution model.   It ensures the existence of the
estimation criterion.

\subsection{Definition of the estimator}
\label{crit}\label{critsimple}

As already mentioned in the introduction,
the starting point of our estimation procedure is to construct an
estimator of the least square contrast
\begin{eqnarray}
\label{Stheta}
 S_{\theta^0,P_X}(\theta)=\mathbb{E}[(Z_1-f_{\theta}(X_{0}))^2 \, w(X_{0})]
,\end{eqnarray}
based on the
observations $(Z_i)$ for $i=0,\ldots,n$.

We  consider the following condition:
 there exists a
weight function $w$ such that for all $\theta\in \Theta$,
\begin{align}
&\condition{C11}\mbox{ The functions }(wf_\theta) \mbox{ and } (wf_\theta^2) \mbox{ belong to
  }\mathbb{L}_1(\mathbb{R}), \mbox{ and  the functions }
w^*/f_\varepsilon^*\,,(f_\theta
w)^*/f_\varepsilon^*\,,\\&\notag
\,(f_\theta^2
w)^*/f_\varepsilon^*  \mbox{ belong to  } \mathbb{L}_1(\mathbb{R}).
\end{align}
\begin{remark}
The first part of Condition  \eref{C11} is not restrictive. The second
part  can be  heuristically expressed  as ``one can  find  a weight
function  $w$ such  that $wf_\theta$  is smooth  enough  compared to
$f_\varepsilon$''.
For a large number of regression functions, such a weight function
can  be easily  exhibited. Some practical choices
are discussed in the
simulation study (Section \ref{simu_context}). 
\end{remark}

If \eref{C11} holds, the expectations
$\mathbb{E}(w(X))$,          $\mathbb{E}(w(X)f_\theta(X))$         and
$\mathbb{E}(w(X)f_\theta^2(X))$  can  be   easily  estimated.  Let  us
present the ideas of the estimation procedure.
 Let $\varphi$ be such that
$\varphi$ and $ \varphi^*/f_\varepsilon^* $ belong to $\mathbb{L}_1(\mathbb{R})$.
For such a function,  due to the independence between $\varepsilon_0$
and $X_0$ we have
\begin{eqnarray*}
\mathbb{E}[\varphi(X_0)]=\mathbb{E}\left( \frac{1}{2\pi}\int \varphi^*(t)e^{-itX_0}dt\right)=\mathbb{E}\left( \frac{1}{2\pi}\int \frac{\varphi^*(t)e^{-itZ_0}}{f_\varepsilon^*(-t)}dt\right).
\end{eqnarray*}
Hence, based on the observations $Z_0,\cdots,Z_n$, $\mathbb{E}[\varphi(X_0)]$ is estimated by
$$\frac{1}{2\pi }\mathbb{R}e\int \frac{\varphi^*(t)n^{-1}\sum_{j=1}^n e^{-itZ_j}}{f_\varepsilon^*(-t)}dt.$$
We  then  propose  to  estimate  $  S_{\theta^0,P_X}(\theta)$  by  the
quantity $S_n(\theta)$ defined
by
\begin{equation}\label{Sn}
 S_{n}(\theta)
=
\frac{1}{2\pi n}\sum_{k=1}^n
\mathbb{R}e\int \frac{\left(\big(Z_k-f_{\theta}\big)^2w\right)^*(t)\; e^{-itZ_{k-1}}}{f_\varepsilon^*(-t)}dt,
\end{equation}
which satisfies
$$\mathbb{E}(S_n(\theta))=\mathbb{E}[(Z_1-f_\theta(X_0))^2w(X_0)].$$
This criteria  is minimum when $\theta=\theta^0$ under the identifiability assumption \eref{Spttilde}.
Using this empirical criterion we propose to estimate $\theta^0$ by
\begin{eqnarray}
\widehat{\theta}=\arg\min_{\theta \in
\Theta}
S_{n}(\theta).\label{thetac}
\end{eqnarray}

\section{Asymptotic properties}
\label{proprietes}

In this section, we give some conditions under which our estimator is consistent and asymptotically normal.

\subsection{Consistency of the estimator}
\label{secparametric}

The first  result to mention is  the consistency of  our estimator. It
holds under the following additional condition.
\begin{align}
&\condition{C11d1} \mbox{ The functions }\sup_{\theta\in \Theta}\left\vert\big( f^{(1)}_{\theta, i} w\big)^*/f_\varepsilon^*\right\vert \mbox{ and }  \sup_{\theta\in
  \Theta}\left\vert\big( f_\theta f^{(1)}_{\theta, i} w\big)^*/f_\varepsilon^* \right\vert\mbox{   belong  to  } \mathbb{L}_1(\mathbb{R})\mbox{ for any } \\&
  \notag
   i \in \{1, \ldots, d\}.
\end{align}
This condition is similar to \eref{C11} for the first
derivatives of $f_\theta$. Thus it is not more restrictive than \eref{C11}.

\begin{theo}\label{consist}
Consider Model \eref{model} under the assumptions
\eref{df}-\eref{emotilde}, \eref{Spttilde},  \eref{emo2tilde},  \eref{fepsnn},  and  the  conditions \eref{C11}-\eref{C11d1}.
Then $\widehat{\theta}$ defined by
\eref{thetac} converges in probability to
$\theta^0$.

\end{theo}

\subsection{$\sqrt{n}$-consistency and asymptotic normality}
To state the asymptotic normality of our estimator,
we need to introduce some  additional  conditions.
\begin{align}
 &\condition{C11d2}  \mbox{  the  functions  } \sup_{\theta\in
  \Theta}\left|\left(
     f^{(2)}_{\theta, i, j} w
 \right)^*/f_\varepsilon^* \right |           \mbox{           and           }
\sup_{\theta\in
  \Theta} \left | \left(\frac{\partial^2}{\partial \theta_i \partial\theta_j} (f_{\theta}^2w)
 \right)^*/f_\varepsilon^* \right| \mbox{        belong       to       }
 \mathbb{L}_1(\mathbb{R}) \mbox{ for}\\&\notag
  \mbox{ any }  i, j \in \{1, \ldots, d\};\\
&\condition{C11d3}  \mbox{  the
  functions                                            }\sup_{\theta\in
  \Theta}\left\vert\left(  \frac{\partial^3
    (f_\theta w)}{\partial \theta_i\partial \theta_j\partial\theta_k}
 \right)^*/f_\varepsilon^*\right\vert  \mbox{  and  }  \sup_{\theta\in
 \Theta}\left\vert\left(\frac{\partial^3}{\partial\theta_i\partial\theta_j\partial\theta_k}(f_\theta ^2 w)
 \right)^*/f_\varepsilon^*\right\vert \mbox{     belong     to     }   \\
&\notag
  \mathbb{L}_1(\mathbb{R}), \mbox{ for } i,j,k \in \{1,\cdots,d\} .\\
&\condition{C11d4}    \mbox{The    integrals     }    \int    \vert    t
(f_{\theta^0}w)^*(t)\vert dt \mbox{ and } \int    \vert    t
(f_{\theta^0}f_{\theta^0,k}^{(1)}w)^*(t)\vert dt \mbox{ are finite,
  for }k\in\{1,\cdots,d\}.
\end{align}

The asymptotic properties of $\widehat{\theta}$, defined by
\eref{thetac}, are stated under two different dependency
conditions, which are presented below.

\begin{defi}\label{defi} Let $(\Omega, {\mathcal A}, {\mathbb P})$ be a probability space.
Let $Y$ be a random variable with values in a Banach space
$({\mathbb B}, \|\cdot \|_{\mathbb B})$. Denote by $\Lambda_\kappa ({\mathbb B})$ the set of
$\kappa$-Lipschitz functions, \textit {i.e.} the  functions $f$ from $({\mathbb B}, \|\cdot\|_{\mathbb B})$ to
${\mathbb R}$ such that $\vert
f(x)-f(y)\vert \leq \kappa \parallel x-y\parallel_{\mathbb{B}}$. Let ${\mathcal M}$ be a
$\sigma$-algebra of ${\mathcal A}$. Let ${\mathbb P}_{Y|{\mathcal
M}}$ be a conditional distribution of $Y$ given ${\mathcal M}$, $\mathbb{P}_Y$
the distribution of $Y$, and  ${\mathcal B}({\mathbb
B})$ the Borel $\sigma$-algebra on $({\mathbb B}, \| \cdot
\|_{\mathbb B})$.
The  dependence coefficients $\alpha$ and $\tau$ are defined by
\begin{eqnarray*}
\alpha({\mathcal M}, \sigma(Y))&=& \frac{1}{2} \sup_{A \in
{\mathcal B}({\mathbb B})}{\mathbb E}(| {\mathbb P}_{Y|{\mathcal
M}}(A)-{\mathbb P}_Y(A)| ) \, , \\
\text{and if ${\mathbb E}(\|Y\|_{\mathbb B}) < \infty$,} \quad \tau ({\mathcal
M}, Y)&=& {\mathbb E}\Big (\sup_{f \in {\Lambda_1}({\mathbb B})}
|{\mathbb P}_{Y|{\mathcal M}}(f)-{\mathbb P}_Y(f)| \Big ) \, .
\end{eqnarray*}

Let ${\bf X}=(X_i)_{i \geq 0}$ be a strictly stationary Markov chain
of real-valued random variables.
On
${\mathbb R}^2$, we put the norm $\|x\|_{{\mathbb
R}^2}=(|x_1|+|x_2|)/2$.
For any integer $k \geq 0$, the coefficients $\alpha_{{\bf
X}}(k)$ and $\tau_{{\bf X}, 2}(k)$ of the chain are defined by
\begin{eqnarray*}
\alpha_{{\bf X}}(k)&=&
 \alpha(\sigma(X_0), \sigma(X_{k}))\\
\text{and if ${\mathbb E}(|X_0|) < \infty$,} \quad \tau_{{\bf X}, 2}(k)&=&  \sup \left \{
 \tau(\sigma(X_0), (X_{i_1},  X_{i_2})), k \leq i_1 \leq  i_2 \right \}. \\
\end{eqnarray*}
\end{defi}
Coefficient $\alpha({\mathcal M}, \sigma(Y))$ is the usual strong  mixing   coefficient
introduced  by Rosenblatt \citeyear{Rosenblatt}. Coefficient $\tau({\mathcal M}, Y)$ has been
introduced by Dedecker and Prieur \citeyear{JDCPCoeff}.
In Section \ref{mixprop}, we recall some
conditions on  $\xi_0$ and $f_{\theta^0}$  under which the  Markov chain
$(X_i)_{i\geq   0}$  is   $\alpha$-mixing   or  $\tau$-dependent   and
illustrate those conditions through some examples.

First we state the asymptotic
normality of $\widehat{\theta}$ when the Markov chain $(X_i)$ of Model \eref{model} is $\alpha$-mixing.
\begin{theo}
\label{NAalpha} Consider Model \eref{model} under assumptions
\eref{df}, \eref{emotilde},
\eref{Spttilde},  \eref{emo2tilde},  \eref{fepsnn},  and conditions  \eref{C11}-\eref{C11d3}.
Let $Q_{|X_1|}$ be the inverse cadlag of the tail function
$t
\to {\mathbb P}(|X_1|>t)$.
Assume  that
\begin{equation}\label{condalpha}
\sum_{k\geq 1}\int_0^{\alpha_{\bf X}(k)}Q_{|X_1|}^2(u) du<\infty \, .
\end{equation}
Then $\widehat{\theta}$ defined by
\eref{thetac} is a $\sqrt{n}$-consistent estimator of
$\theta^0$ which satisfies
$$\sqrt{n}(\widehat{\theta}-\theta^0)\cvl              \mathcal{N}(0,
\Sigma_1),$$
where the covariance matrix $\Sigma_1$  is defined in equation \eref{secondequality}.
\end{theo}

Next, we give the corresponding result when the Markov chain $(X_i)$ is $\tau$-dependent.

\begin{theo}
\label{NAtau}
Consider Model \eref{model} under assumptions
\eref{df}, \eref{emotilde},
\eref{Spttilde},  \eref{emo2tilde},  \eref{fepsnn},  and conditions  \eref{C11}-\eref{C11d4}.
Let $G(t)=t^{-1}{\mathbb E}(X_1^2{\mathbf 1}_{X_1^2>t})$,  and let $G^{-1}$ be the inverse cadlag of $G$.
Assume  that
\begin{equation}\label{condtau}
\quad \sum_{k>0}G^{-1}(\tau_{{\bf X}, 2}(k)) \tau_{{\bf X}, 2}(k) < \infty  \, .
\end{equation}
Then $\widehat{\theta}$ defined by
\eref{thetac} is a $\sqrt{n}$-consistent estimator of
$\theta^0$ which satisfies
$$\sqrt{n}(\widehat{\theta}-\theta^0)\cvl              \mathcal{N}(0,
\Sigma_1), $$
where the covariance matrix $\Sigma_1$  is defined in equation \eref{secondequality}.
\end{theo}

\begin{remark} Let us give some  conditions under which (\ref{condalpha}) or (\ref{condtau}) are verified. Assume that ${\mathbb E}(|X_0|^p)< \infty$  for some $p>2$.
Then  (\ref{condalpha}) is true provided that $\sum_{k>0} k^{2/(p-2)} \alpha_{\bf X}(k)< \infty$, and  (\ref{condtau}) is true provided that
$\sum_{k>0} (\tau_{{\bf X}, 2}(k))^{(p-2)/p}<\infty$.
\end{remark}

Note  that  those  results do  not  require  the  Markov chain  to  be
absolutely  regular   as  it  is   the  case  in  Comte   and  Taupin
\citeyear{ComteTaupin2001}. Consequently they apply to autoregressive
models  with  weaker  dependency  conditions.  Beside  the  dependency
conditions, our estimation procedure allows to achieve
the parametric rate for a larger class of regression functions than  in  Comte   and  Taupin
\citeyear{ComteTaupin2001}.

The conditions under which Theorems \ref{NAalpha} and \ref{NAtau} hold
are  similar,  except Condition  \eref{C11d4}  which appears only in Theorem \ref{NAtau}. This condition is  just
technical and not restrictive at all.

The choice of the weight function $w$ is crucial.
Various
weight functions can handle with Conditions \ref{C11}-\ref{C11d4}.  The numerical
properties of  the resulting estimators  will differ from one  choice to
another.  This  point  is  discussed  on  simulated data  in  the next section.

\section{Simulation study}
\label{simu_context}
We investigate the properties of our estimator for different regression functions on simulated data. For each choice of regression function, we consider
two error distributions:
 the
Laplace distribution
 and the
Gaussian distribution. 
When  $\varepsilon_1$  has  the  Laplace distribution,  its  density  and
Fourier transform are
\begin{equation}\label{doublexp}
f_{\varepsilon}(x)= \frac{1}{\sigma_\varepsilon\sqrt{2}}
\exp\Big(-\frac{\sqrt{2}}{\sigma_\varepsilon}|x|\Big), \mbox { and } f^*_{\varepsilon}(x)= \frac{1}{1+\sigma_\varepsilon^2 x^2/2}.
\end{equation}
Hence, $\varepsilon_1$ is centered
with variance $\sigma_\varepsilon^2$.

\medskip

\noindent  When $\varepsilon_1$  is  Gaussian, its  density and  Fourier
transform are
\begin{equation}\label{gauss}
f_{\varepsilon}(x)= \frac{1}{\sigma_\varepsilon\sqrt{2\pi}}
\exp\Big(-\frac{x^2}{2\sigma_\varepsilon^2}\Big), \mbox { and } f^*_{\varepsilon}(x)=
\exp(-\sigma_\varepsilon^2 x^2/2).
\end{equation}
Hence,  $\varepsilon_1$ is centered
with variance $\sigma_\varepsilon^2$.

\medskip
For each of these error distributions, we consider the case of a linear regression function and of a Cauchy regression function. We start with the linear case.

\subsection{Linear regression function}
\label{regsetting}
We consider the model (\ref{model})  with $f_\theta(x)= ax + b$, where $|a|<1$ and
$\theta=(a,b)^T$.
In these simulations, we have  chosen to illustrate the numerical properties
of  our estimator  under the  weakest of the dependency  conditions,  that is
$\tau$-dependency.
As it is recalled in Appendix \ref{mixprop}, when $f_{\theta^0}$ is linear with $|a|<1$, if
$\xi_0$ has a density bounded from below in a neighborhood of the origin, then  the
Markov chain $(X_i)_{i \geq 0}$ is $\alpha$-mixing. When   $\xi_0$  does not  have  a  density, then  the  chain  may not be
$\alpha$-mixing (and not even  irreducible), but it is always $\tau$-dependent.

Here, we consider the case where the    innovation
distribution is discrete, in such a way that the stationary Markov Chain  is $\tau$-dependent but not $\alpha$-mixing.
We also consider two distinct values of  $\theta_0$.
For the first value, the
stationary  distribution of $X_i$ is absolutely continuous with respect to the Lebesgue measure. For the second value,
the stationary distribution is singular with respect to the Lebesgue measure.
In both cases Theorem \ref{NAtau} applies, and the estimator $\hat \theta$ is asymptotically normal.

\medskip

\noindent $\bullet$ {Case A (absolutely continuous stationary distribution).}
We focus on the case where the
true parameter is $\theta^0=(1/2, 1/4)^T$,  $X_0$ is uniformly distributed over $[0, 1]$, and  $(\xi_i)_{i \geq 1}$ is a sequence of i.i.d. random variables, independent of $X_0$
and such that ${\mathbb P}(\xi_1=-1/4)={\mathbb P}(\xi_1=1/4)=1/2$. Then the Markov chain defined for $i>0$ by
\begin{equation}
\label{lin1}   X_i= \frac14+\frac12 X_{i-1} + \xi_i
\end{equation}
is strictly stationary, the stationary distribution
being the uniform distribution over $[0, 1]$, and consequently $\sigma_{X_0}^2=1/12$.
 This chain  is non-irreducible,  and the dependency  coefficients are
 such  that  $\alpha_{\bf  X}  (k)=  1/4$ (see  for  instance  Bradley
 \citeyear{Bradley86},     p.     180)     and     $\tau_{{\bf     X},
   2}(k)=O(2^{-k})$. Thus the Markov  chain is not $\alpha$-mixing, but it
 is $\tau$-dependent.
For the simulation, we
start  with $X_0$  uniformly distributed  over $[0,  1]$, so  the
simulated chain is stationary.

\medskip

\noindent $\bullet$ {Case B (singular stationary distribution).}
We consider the case where the
true parameter is $\theta^0=(1/3, 1/3)^T$,  $X_0$ is uniformly distributed over the Cantor set,  and  $(\xi_i)_{i \geq 1}$ is a sequence of i.i.d. random variables, independent of $X_0$
and such that ${\mathbb P}(\xi_1=-1/3)={\mathbb P}(\xi_1=1/3)=1/2$. Then the Markov chain defined for $i>0$ by
\begin{equation}
\label{lin2}   X_i= \frac13+\frac13 X_{i-1} + \xi_i
\end{equation}
is strictly stationary, the stationary distribution
being the uniform distribution over the Cantor set, and consequently  $\sigma_{X}^2=1/8$. This chain is non-irreducible, and the dependency coefficients
satisfy    $\alpha_{\bf   X}   (k)=    1/4$   and    $\tau_{{\bf   X},
  2}(k)=O(3^{-k})$. Thus the Markov  chain is not $\alpha$-mixing, but
 is $\tau$-dependent.
For the simulation, we start with $X_0$ uniformly distributed over $[0, 1]$, and we consider that the chain is close
to  the   stationary  chain  after   1000  iterations.  We   then  set
$X_i=X_{i+1000}$.

\medskip

In these two cases, we can  find a weight function $w$ satisfying the conditions \eref{C11}-\eref{C11d4}.
We first give the detailed expression  of the estimator for two choices of  weight functions $w$. Then we recall the classic estimator when $X$ is directly observed, the ARMA estimator, and the so-called naive estimator.

\subsubsection{Expression of the estimator.}
\label{Linear}
We consider the  two following weight
functions $w$
\begin{eqnarray}
\label{weightlin}
w(x)=N(x)=\exp\{-x^2/(4 \sigma_\varepsilon^2)\} \ \mbox{ and } \ w(x)=SC(x)=\frac{1}{2\pi}\Big(\frac{2*\sin(x)}{x}\Big)^4.\end{eqnarray}
 These  choices of weight  ensure that  Conditions \eref{C11}-\eref{C11d4}
 hold and that  the two estimators, denoted by $\widehat{\theta}_N$ and $\widehat{\theta}_{SC}$ respectively, converge  to $\theta^0$ with  the parametric
 rate of convergence. There are two main differences between these two weight
 functions.    First,   $N$    depends   on    the    variance   error
 $\sigma_\varepsilon^2$. Hence the estimator should be adaptive to the
 noise  level. On  the contrary,  it may  be sensitive  to  very small
 error  variance  as  it  appears  in  the  simulations  (see  Figure
 \ref{fig0}).
Second,  $SC$  has  strong  smoothness properties  since  its  Fourier
transform is compactly supported.

The    two associated  estimators  are  based  on the calculation of  $S_n(\theta)$, which can be
written as
$$S_n(\theta)=\frac{1}{n}\sum_{k=1}^n    [(Z_k^2+b^2-2Z_kb)I_0(Z_{k-1})
+a^2I_2(Z_{k-1})-2a(Z_k-b)I_1(Z_{k-1})],$$
with
\begin{equation}\label{Ij}
I_j(Z)=\frac{1}{2\pi}\mathbb{R}e\int      (p_jw)^*(u)\frac{e^{-iuZ}}{f_\varepsilon^*(-u)}du,
\end{equation}
where $p_j(x)=x^j$ for $j=0,1,2$,
$w$ being either $w=N$ or $w=SC$.
With the above notations,
 $\widehat{\theta} = (\widehat{a}, \widehat{b})^T$
satisfies
\begin{eqnarray}
\label{coeflin}
\widehat{a}&=&
\frac{\sum_{k=1}^n Z_k I_{1}(Z_{k-1})\sum_{k=1}^nI_{0}(Z_{k-1})-
\sum_{k=1}^n Z_k I_{0}(Z_{k-1})\sum_{k=1}^nI_{1}(Z_{k-1})
}{\sum_{k=1}^n I_{2}(Z_{k-1})\sum_{k=1}^n I_{0}(Z_{k-1})-\big(\sum_{k=1}^n I_{1}(Z_{k-1})\big)^2},\\
\label{coeflin2}
\widehat{b}&=&\frac{\sum_{k=1}^n Z_k I_{0}(Z_{k-1})}{\sum_{k=1}^nI_{0}(Z_{k-1})}-\widehat{a}\frac{
\sum_{k=1}^n I_{1}(Z_{k-1})}{\sum_{k=1}^nI_{0}(Z_{k-1})}.
\end{eqnarray}
We now compute $I_j(Z)$ for $j=0,1,2$ and the two weight functions.
In the following we respectively denote $I_{j,N}(Z)$ and $I_{j,SC}(Z)$ the previous
integrals when the weight function is either $w=N$ or $w=SC$.

We start with $w = N$ and give the details of the calculations for the two error distributions (Laplace and Gaussian), which are explicit.  Then, with the weight function $w = SC$, we present the calculations, which are not explicit whatever the error distribution $f_\varepsilon$.

\medskip

\noindent $\bullet$ When $w=N$,
Fourier calculations provide that
\begin{eqnarray*}
N^*(t)&=&\sqrt{2\pi}\sqrt{2\sigma_\varepsilon^2}\exp(-\sigma_\varepsilon^2 t^2)\\
 (N p_1)^*(t)&=&\sqrt{2\pi}\sqrt{2\sigma_\varepsilon^2}\exp(-\sigma_\varepsilon^2
 t^2) \big(-2\sigma_\varepsilon^2 t/i\big),\\
 (N p_2)^*(t)&=&-\sqrt{2\pi}\sqrt{2\sigma_\varepsilon^2}\exp(-\sigma_\varepsilon^2
 t^2)\big(-2\sigma_\varepsilon^2+4\sigma_\varepsilon^4t^2\big).
\end{eqnarray*}

\noindent It follows that
\begin{eqnarray*}
I_{0,N}(Z)&=& \frac{1}{2\pi}\mathbb{R}e\int  \sqrt{2\pi}\sqrt{2\sigma_\varepsilon^2}\exp(-\sigma_\varepsilon^2 t^2)\frac{e^{-itZ}}{f_\varepsilon^*(-t)}dt,\\
I_{1,N}(Z)&=& \frac{1}{2\pi}\mathbb{R}e\int \sqrt{2\pi}\sqrt{2\sigma_\varepsilon^2}\exp(-\sigma_\varepsilon^2
 t^2) \big(-2\sigma_\varepsilon^2 t/i\big)\frac{e^{-itZ}}{f_\varepsilon^*(-t)}dt,
\\
I_{2,N}(Z)&=& \frac{1}{2\pi} \mathbb{R}e\int \sqrt{2\pi}\sqrt{2\sigma_\varepsilon^2}\exp(-\sigma_\varepsilon^2
 t^2)\big(2\sigma_\varepsilon^2-4\sigma_\varepsilon^4t^2\big)\frac{e^{-itZ}}{f_\varepsilon^*(-t)}dt.
 \end{eqnarray*}

\noindent If $f_\varepsilon$ is the Laplace distribution (\ref{doublexp}), replacing
$f_\varepsilon^*$ by its expression
we get
\begin{eqnarray*}
I_{0,N}(Z)
&=&  e^{-Z^2/(4\sigma_\varepsilon^2)}-\frac{\sigma_\varepsilon^2}{2}\frac{\partial^2 }{\partial Z^2}N(Z)
=  \left[5/4-Z^2/(8\sigma_\varepsilon^2)\right]e^{-Z^2/(4\sigma_\varepsilon^2)},
\end{eqnarray*}
\begin{eqnarray*}
I_{1,N}(Z)=\left[7Z/4-Z^3/(8\sigma_\varepsilon^2)\right]e^{-Z^2/(4\sigma_\varepsilon^2)},
I_{2,N}(Z)=\left[-\sigma_\varepsilon^2+9Z^2/4-Z^4/(8\sigma_\varepsilon^2) \right]e^{-Z^2/(4\sigma_\varepsilon^2)}.\end{eqnarray*}

\noindent If $f_\varepsilon$ is the Gaussian distribution (\ref{gauss}), replacing
$f_\varepsilon^*$ by its expression we obtain
$$
I_{0,N}(Z)= \sqrt{2}e^{-Z^2/(2\sigma_\varepsilon^2)}, \quad I_{1,N}(Z)=2\sqrt{2}Ze^{-Z^2/(2\sigma_\varepsilon^2)}
\quad \mbox{and} \quad ~~~I_{2,N}(Z)=\sqrt{2}(4Z^2-2\sigma_\varepsilon^2) e^{-Z^2/(2\sigma_\varepsilon^2)}.
$$
Hence we deduce the  expression of $\widehat{a}_N$ and $\widehat{b}_N$
by applying \eref{coeflin} and \eref{coeflin2}.

\medskip
\noindent $\bullet${When $w=SC$}, Fourier calculations provide that
\begin{eqnarray*}
SC^*(t)&=&\ind_{[-4,-2]}(t)(t^3/6+2t^2+8t+32/3)+\ind_{[-2,0]}(t)(-t^3/2-2t^2+16/3)\\
&&
+\ind_{[2,4]}(t)(-t^3/6+2t^2-8t+32/3)+\ind_{[0,2]}(t)(t^3/2-2t^2+16/3)\\
 (SC p_1)^*(t)&=&\frac{\partial}{\partial t}SC^*(t)/i \mbox{ and }
 (SC p_2)^*(t)=\frac{\partial^2}{\partial t^2}SC^*(t)/(i^2).
\end{eqnarray*}
The integrals $I_{j,SC}(Z)$, defined for $j=0,1,2$ by
\begin{eqnarray}
\label{ISClin}
I_{j,SC}(Z)=\frac{1}{2\pi}\mathbb{R}e\int
(SC p_j)^*(t)\frac{e^{-itZ}}{f_\varepsilon^*(-t)}dt,
\end{eqnarray}
  have no  explicit  form, whatever the error distribution $f_\varepsilon$.  It has  to be  numerically
computed, using the  IFFT Matlab function.  More precisely, we  consider a
finite         Fourier          series         approximation         of
${(SC p_j)^*(t)}/{f_\varepsilon^*(t)}$ whose Fourier transfom is
calculated  using IFFT  Matlab function.  The  result is  taken as  an
approximation of $I_{j,SC}(Z)$. Finally we deduce the  expression of $\widehat{a}_{SC}$ and $\widehat{b}_{SC}$
by applying \eref{coeflin}  and \eref{coeflin2}.

\subsubsection{Comparison with classical estimators}
We compare the two estimators $\widehat{\theta}_{N}$ and $\widehat{\theta}_{SC}$ with three classical estimators, the usual least square  estimator when there is no observation noise, the ARMA estimator, and the so-called naive estimator.

\medskip

\noindent $\bullet$  {Estimator without noise.}
In the case where $\varepsilon_i=0$, that is $(X_0, \ldots, X_n)$ is observed without error, the parameters can be easily estimated
by the usual least square estimators
$$\widehat{a}_{X}=\frac{n\sum_{i=1}^n X_i X_{i-1} - \sum_{i=1}^n X_i \sum_{i=1}^n X_{i-1}}
{n\sum_{i=1}^n X_{i-1}^2- (\sum_{i=1}^n X_{i-1})^2} \quad  \mbox{and} \quad \widehat{b}_{X}=\frac 1n \Big(\sum_{i=1}^n X_i\Big)-
\widehat{a}_{X}\frac1n\Big(\sum_{i=1}^n X_{i-1}\Big)\, .$$

\medskip

\noindent $\bullet$ {ARMA estimator.}
When the regression function is linear, the model may be written as
$$
Z_i- a Z_{i-1}-b= \xi_i + \varepsilon_i -a \varepsilon_{i-1} \, .
$$
The auto-covariance function $\gamma_Y$ of the stationary sequence
$Y_i=\xi_i + \varepsilon_i -a \varepsilon_{i-1}$ is given
by
$$
  \gamma_Y(0)=(1+a^2)\sigma_\varepsilon^2+ \sigma_\xi^2, \quad \gamma_Y(1)= -a \sigma_\varepsilon^2, \quad \text{and $\gamma_Y(k)=0$ for $k>1$}.
$$
It follows that $Y_i$ is an MA(1) process, which may be written as
$$
  Y_i= \eta_i-\beta\eta_{i-1},
$$
where $\eta_i$ is the innovation, and $|\beta|<1$ (note that $|\beta|\neq 1$
because $\gamma_Y(0)-2|\gamma_Y(1)| > 0$). Moreover, one can give the explicit
expression of
$\beta$ and  $\sigma_\eta^2$  in terms of $a, \sigma_\xi^2$ and $\sigma_\varepsilon^2$.
It follows that, if $|a|<1$, $(Z_i)_{i \geq 0}$ is the causal invertible
ARMA(1,1) process
\begin{eqnarray}
\label{arma}
Z_i-aZ_{i-1}=b+\eta_i-\beta\eta_{i-1}.
\end{eqnarray}

Note that $a\neq \beta$ except if $a=0$. Hence, if $|a|<1$ and $a\neq 0$,
one  can  estimate the parameters $(a, b, \beta)$  by  maximizing the  so-called
Gaussian likelihood.
These estimators are
consistent and asymptotically Gaussian. Moreover they are
efficient when both the innovations and the errors
$\varepsilon$   are  Gaussian   (see  Hannan   \citeyear{Hannan73} or
Brockwell and Davis \citeyear{BrockwellDavis}).
Note that this well-known approach
does not require the knowledge of the error
distribution, but of course it works only in the particular case
where the regression function $f_\theta$ is
linear.  For the computation of the ARMA estimator we use the function \textit{arma} from the R
\textit{tseries} package (see Trapletti and Hornik \citeyear{Rtseries}).
The resulting
estimators    are      denoted     by     $\widehat{a}_{arma}$     and
$\widehat{b}_{arma}$.

\medskip

\noindent $\bullet$ {Naive estimator.}
 The naive estimator is constructed by replacing the unobserved $X_i$ by the observation $Z_i$
in the expression of $\widehat{a}_X$ and $\widehat{b}_X$:
$$\widehat{a}_{naive}=\frac{n\sum_{i=1}^n Z_i Z_{i-1} - \sum_{i=1}^n Z_i \sum_{i=1}^n Z_{i-1}}
{n\sum_{i=1}^n Z_{i-1}^2- (\sum_{i=1}^n Z_{i-1})^2} \quad  \mbox{and} \quad \widehat{b}_{naive}=\frac 1n \Big(\sum_{i=1}^n Z_i\Big)-
\widehat{a}_{naive}\frac1n\Big(\sum_{i=1}^n Z_{i-1}\Big)\, .$$
Classical results show that $\widehat{\theta}_{naive}$ is an asymptotically  biased estimator of $\theta^0$, which is confirmed by the simulation study.

\subsubsection{Simulation results}
For each error distribution, we  simulate 100   samples with  size  $n$, $n=500$,
$5000$ and $10000$. We consider different values of  $\sigma_\varepsilon$  such that the ratio signal
to noise $s2n=\sigma_\varepsilon^2/\mbox{Var}(X)$
is $0.5,1.5$ or $3$.
The comparison of the five estimators  is  based on the bias, the Mean Squared
Error (MSE), and the box plots.
If $\widehat{\theta}(k)$  denotes the value of the  estimation for the
$k$-th sample, the MSE is evaluated by the empirical mean over the 100 samples:
\begin{eqnarray*}
MSE(\widehat{\theta})=\frac{1}{100}\sum_{k=1}^{100}(\widehat{\theta}(k)-\theta^0)^2.
\end{eqnarray*}
Results are presented in  Figures  \ref{fig0}-\ref{fig00}   and  Tables
\ref{table1}-\ref{table4}.

{\tiny
\begin{table}
\begin{tabular}{lllrrrrr}
\hline
& ratio&& \multicolumn{5}{c}{Estimator}\\
n&s2n&& \multicolumn{1}{c}{$\widehat{\theta}_{arma}(MSE)$} & \multicolumn{1}{c}{$\widehat{\theta}_N (MSE)$} & \multicolumn{1}{c}{$\widehat{\theta}_{SC}(MSE)$}& \multicolumn{1}{c}{$\widehat{\theta}_X (MSE)$}  &\multicolumn{1}{c}{$\widehat{\theta}_{naive}(MSE)$}\\
\hline

$1000$ & 0.5&
 a     & 0.487   (0.008) &   0.459     (0.020) &   0.489    (0.002)&    0.493    (0.001)&    0.328    (0.030)\\
&&b  &    0.257    (0.002)   & 0.262    (0.002)   & 0.255   (0.001) &   0.253    (0.001)  &  0.336    (0.008)\\
\\
&  1.5&
 a   &   0.494    (0.015)  &  0.488    (0.013) &   0.492   (0.006) &   0.501    (0.001)  &  0.198   (0.092)\\
&&b  &    0.251    (0.004) &   0.253    (0.002)  &  0.253    (0.002) &   0.249    (0.001)  &  0.399    (0.023)\\
\\
 &  3&
  a   &   0.461    (0.044)   & 0.502    (0.029)  &  0.503    (0.026)  &  0.493    (0.001)   & 0.121    (0.145)\\
&&b    &  0.270    (0.012) &   0.249   (0.001)  &  0.249    (0.001) &   0.253   (0.001) &   0.440    (0.037)\\
 \hline
$5000$ &  0.5&
 a    &  0.497    (0.001)  &  0.499  (0.004) &   0.499    (0.001) &   0.499   (0.001)  & 0.332    (0.028)\\
&&b &     0.252    (0.001)   & 0.251    (0.001)  &  0.251    (0.001) &   0.251    (0.001) &   0.334    (0.007)\\
\\
&  1.5&
 a    &  0.498    (0.003)  &  0.508  (0.003)  &  0.503    (0.002)  &  0.499    (0.001) &   0.199    (0.091)\\
&&b   &   0.250    (0.001)  &  0.247    (0.001)  &  0.248    (0.001)  &  0.250    (0.001)  &  0.399    (0.022)\\
\\
&  3&
 a   &   0.487    (0.008) &   0.492    (0.004)  &  0.495    (0.004)&    0.500    (0.001)  &  0.123    (0.143)\\
&&b  &    0.256    (0.002) &   0.253    (0.001) &  0.252    (0.001)  &  0.250    (0.001)    &0.437   (0.035)\\
\hline
$10000$ & 0.5&
 a     & 0.496    (0.001)  &  0.501    (0.002) &   0.500    (0.001)   & 0.499    (0.001)   & 0.334    (0.028)\\
&&b    &  0.252    (0.001)   & 0.250    (0.001)   & 0.250    (0.001)  &  0.250    (0.001)  &  0.333    (0.007)\\
 \\
 & 1.5&
  a  &    0.504    (0.002)  &  0.500    (0.001) &   0.501    (0.001)  &  0.500    (0.001)  & 0.200   (0.090)\\
&&b   &   0.248    (0.001) &   0.250   (0.001)   & 0.250    (0.001)  &  0.250    (0.001) &   0.401   (0.023)\\
   \\
 & 3&
  a    &  0.493    (0.003)   & 0.499    (0.001)   & 0.499    (0.002) &   0.498    (0.001)  &  0.124    (0.142)\\
&&b   &   0.254    (0.001)  &  0.250    (0.001)  &  0.250    (0.001)   & 0.251    (0.001)   & 0.438    (0.036)\\

 \hline\\
\end{tabular}

\caption{\small{Estimation results for Linear Case  A, Laplace error. Mean estimated values of the five estimators $\widehat{\theta}_{arma}$, $\widehat{\theta}_{N}$, $\widehat{\theta}_{SC}$, $\widehat{\theta}_{X}$ and $\widehat{\theta}_{naive}$ are presented for various values of  $n$ ($1000$,  5000 or 10000) and s2n (0.5, 1.5, 3). True values are $a^0
  = 1/2$, $b^0 = 1/4$. MSEs are given in brackets.}}
\label{table1}
\end{table}}

{\tiny
\begin{table}
\begin{tabular}{lllrrrrr}
\hline
& ratio&& \multicolumn{5}{c}{Estimator}\\
n&s2n&& \multicolumn{1}{c}{$\widehat{\theta}_{arma}(MSE)$} & \multicolumn{1}{c}{$\widehat{\theta}_N (MSE)$} & \multicolumn{1}{c}{$\widehat{\theta}_{SC}(MSE)$}& \multicolumn{1}{c}{$\widehat{\theta}_X (MSE)$}  &\multicolumn{1}{c}{$\widehat{\theta}_{naive}(MSE)$}\\
\hline

$1000$ & 0.5
 &a&      0.483    (0.006)&    0.539    (0.039)&    0.496     (0.002)&    0.495    (0.001)&    0.331    (0.030)\\
&&b &     0.259     (0.002)&    0.243    (0.003)&    0.253    (0.001)&    0.253    (0.001)&    0.336    (0.008)\\
\\
&  1.5&
 a    &  0.497    (0.021)&    0.516    (0.027)&    0.507     (0.009)&    0.499     (0.001)&    0.200     (0.091)\\
&&b &     0.251     (0.005)&    0.243    (0.005)&    0.246    (0.002)&    0.249     (0.001)&    0.399     (0.023)\\
\\
 &  3&
  a   &   0.456    (0.031)&    0.521    (0.082)&    0.481    (0.030)&    0.501     (0.001)&    0.120    (0.145)\\
&&b &     0.272     (0.008)&    0.244    (0.016)&    0.260     (0.007)&    0.250    (0.001)&    0.441    (0.037)\\
\hline
$5000$ &  0.5&
 a      &0.497    (0.001)&    0.492     (0.006)&    0.499    (0.001)&    0.498     (0.001)&    0.333    (0.028)\\
&&b   &   0.251     (0.001)&    0.252     (0.001)&   0.250     (0.001)&    0.250     (0.001)&    0.333     (0.007)\\
\\
&  1.5&
 a   &   0.490    (0.002)&    0.510    (0.006)&    0.502    (0.001)&    0.499     (0.001)&    0.120    (0.090)\\
&&b  &    0.254     (0.001)&    0.245     (0.001)&    0.248     (0.001)&    0.250    (0.001)&   0.399    (0.022)\\
\\
&  3&
 a    &  0.471     (0.010)&    0.512    (0.008)&    0.503    (0.005)&    0.498     (0.001)&    0.124     (0.141)\\
&&b &     0.263     (0.002)&    0.245     (0.002)&    0.249    (0.001)&    0.251    (0.001)&    0.437    (0.035)\\
\hline
$10000$ & 0.5&
 a   &   0.504   (0.006)&    0.500    (0.003)&    0.498    (0.001)&    0.499    (0.001)&    0.331     (0.028)\\
&&b  &    0.249    (0.001)&   0.250     (0.001)&    0.251     (0.001)&   0.251    (0.001)&    0.335    (0.007)\\
 \\
 & 1.5&
  a    &  0.495     (0.002)&    0.501    (0.002)&    0.499     (0.001)&    0.501      (0.001)&    0.200    (0.090)\\
&&b &     0.253    (0.001)&    0.250    (0.001)&    0.251    (0.001)&    0.250   (0.001)&    0.401    (0.023)\\
  \\
 & 3&
  a  &    0.492     (0.004)&   0.498    (0.004)&    0.500    (0.003)&    0.500    (0.001)&    0.126    (0.140)\\
&&b &     0.254    (0.001)&    0.251    (0.001)&    0.251    (0.001)&    0.250    (0.001)&    0.437    (0.009)
 \\
 \hline\\

\end{tabular}
\caption{\small{Estimation results for Linear Case A, Gaussian error. Mean estimated values of the five estimators $\widehat{\theta}_{arma}$, $\widehat{\theta}_{N}$, $\widehat{\theta}_{SC}$, $\widehat{\theta}_{X}$ and $\widehat{\theta}_{naive}$ are presented for various values of $n$ ($1000$, 5000 or 10000) and s2n (0.5, 1.5, 3). True values are $a^0
  = 1/2$, $b^0 = 1/4$. MSEs are given in brackets.}}
\label{table2}
\end{table}
}

{\tiny

\begin{table}
\begin{tabular}{lllrrrrr}
\hline
& ratio&& \multicolumn{5}{c}{Estimator}\\
n&s2n&& \multicolumn{1}{c}{$\widehat{\theta}_{arma}(MSE)$} & \multicolumn{1}{c}{$\widehat{\theta}_N(MSE)$} & \multicolumn{1}{c}{$\widehat{\theta}_{SC}(MSE)$}& \multicolumn{1}{c}{$\widehat{\theta}_X(MSE)$}  &\multicolumn{1}{c}{$\widehat{\theta}_{naive}(MSE)$}\\
\hline
$1000$ & 0.5
&a&      0.288    (0.021)&  0.341      (0.013)&    0.330     (0.002)&   0.326    (0.001)&  0.217      (0.015)\\
&&b&      0.354    (0.005)&     0.331   (0.001)&    0.333     (0.001)& 0.335      (0.001)&     0.389  (0.004)
 \\
 \\
&  1.5&
a&      0.298    (0.050)&    0.332   (0.009)& 0.335       (0.007)& 0.330       (0.001)&  0.136      (0.040)\\
 &&b&      0.349    (0.012)&  0.331      (0.002)& 0.329      (0.002)&    0.335    (0.001)&    0.429    (0.010) \\
\\
 &  3&
a&      0.240    (0.127)&  0.343      (0.017)&  0.343     (0.018)&   0.330    (0.001)& 0.084       (0.063)\\
  &&b&      0.385     (0.033)&    0.333   (0.003)&    0.333  (0.003)&       0.338  (0.001)&    0.465    (0.018)
 \\
\hline
$5000$ &  0.5&
a&      0.333    (0.004)&    0.335     (0.003)& 0.335       (0.001)&  0.333     (0.001)& 0.223       (0.012)\\
&&b&      0.333     (0.001)& 0.332       (0.001)&   0.332    (0.001)&    0.334    (0.001)&    0.388    (0.003) \\
\\
&  1.5&
a&      0.331    (0.011)&  0.328      (0.002)&    0.334   (0.001)&    0.334    (0.001)&    0.433    (0.041)\\
&&b&    0.334    (0.003)&   0.334     (0.001)&    0.329     (0.001)&    0.332    (0.001)&    0.132    (0.010) \\
\\
&  3&
a&      0.290     (0.030)&  0.329     (0.003)&   0.329     (0.004)&    0.333     (0.001)&  0.083     (0.063)\\
 &&b&      0.355   (0.008)&    0.335    (0.008)&     0.335   (0.008)&    0.334    (0.001)& 0.459     (0.016)
 \\
 \hline
$10000$ & 0.5&
a&      0.337    (0.002)&   0.335     (0.002)&    0.334     (0.001)&    0.334     (0.001)& 0.222      (0.012)\\
 &&b&     0.331     (0.001)&   0.332     (0.001)&    0.332    (0.001)&    0.332     (0.001)&   0.388    (0.003) \\
 \\
 & 1.5&
a&      0.322     (0.006)&  0.336      (0.001)&   0.336    (0.001)&    0.334     (0.001)& 0.134    (0.040)\\
&&b&      0.339    (0.002)&     0.332  (0.001)&    0.332    (0.001)&    0.333     (0.001)&    0.433    (0.010) \\
  \\
 & 3&
a&      0.329    (0.010)&    0.336    (0.002)&    0.336    (0.002)&    0.334     (0.001)&  0.083     (0.063)\\
 &&b&      0.335    (0.002)&    0.332     (0.001)&    0.332     (0.001)&    0.332     (0.001)&  0.457      (0.015)
 \\
 \hline\\

\end{tabular}
\caption{\small{Estimation results for Linear Case B, Laplace error. Mean estimated values of the five estimators $\widehat{\theta}_{arma}$, $\widehat{\theta}_{N}$, $\widehat{\theta}_{SC}$, $\widehat{\theta}_{X}$ and $\widehat{\theta}_{naive}$ are presented for various values of $n$ ($1000$, 5000 or 10000) and s2n (0.5, 1.5, 3). True values are $a^0
  = 1/3$, $b^0 = 1/3$. MSEs are given in brackets.
}}
\label{table3}
\end{table}}

{\tiny
\begin{table}
\begin{tabular}{lllrrrrr}
\hline
& ratio&& \multicolumn{5}{c}{Estimator}\\
n&s2n&& \multicolumn{1}{c}{$\widehat{\theta}_{arma}(MSE)$} & \multicolumn{1}{c}{$\widehat{\theta}_N (MSE)$} & \multicolumn{1}{c}{$\widehat{\theta}_{SC}(MSE)$}& \multicolumn{1}{c}{$\widehat{\theta}_X (MSE)$}  &\multicolumn{1}{c}{$\widehat{\theta}_{naive}(MSE)$}\\
\hline

$1000$ & 0.5
&a  &    0.327    (0.016)&    0.349    (0.035)&   0.330    (0.003) &    0.326    (0.001) &   0.218    (0.014)\\
&& b  &    0.338   (0.004) &   0.332   (0.002) &   0.336    (0.001) &     0.337  (0.001) &  0.392     (0.004)\\
\\
&  1.5&a&    0.290     (0.061) &    0.355   (0.021)&    0.345   (0.008) &   0.332    (0.001) &    0.133   (0.041)\\
 &&b&      0.353     (0.015)&    0.324     (0.004)&     0.328   (0.002)&     0.333   (0.001)&     0.432  (0.010)\\
\\
 &  3
 &a  &    0.234     (0.153)&   0.329     (0.049)&  0.329      (0.051)&   0.326    (0.001)&   0.077    (0.067)\\
&& b &     0.383    (0.040)&   0.337    (0.010)&     0.337  (0.010)&   0.337     (0.001)&    0.461  (0.017)\\
\hline
$5000$ &  0.5&a&      0.329     (0.004)&   0.341    (0.005)&   0.333    (0.001)&     0.332   (0.001)&     0.220  (0.013)\\
 &&b&      0.335     (0.001)&    0.332    (0.001)&  0.334     (0.001)&   0.334    (0.001)&   0.399     (0.003)\\
\\
&  1.5
&a&      0.329     (0.009)&    0.331   (0.003)&    0.332    (0.002)&    0.333     (0.001)&   0.132      (0.041)\\
&&b&      0.335   (0.002)&   0.334     (0.001)&   0.333     (0.001)&    0.333     (0.001)&    0.433    (0.010)\\
\\
&  3
&a&      0.315     (0.022)&    0.348    (0.008)&    0.348     (0.008)&    0.334     (0.001)&   0.084    (0.062)\\
&&b&      0.343    (0.006)&   0.327     (0.002)&    0.328    (0.002)&    0.332     (0.001)&     0.459  (0.016)\\
\hline
$10000$ & 0.5&

a&      0.330    (0.002)&    0.333    (0.003)&    0.333     (0.001)&    0.332     (0.001)&  0.221     (0.013)\\
 &&b&      0.335    (0.001)&    0.333     (0.001)&    0.333    (0.001)&    0.334  (0.001)&   0.389     (0.003)\\
 \\
 & 1.5&
a&      0.328     (0.006)&    0.336     (0.002)&    0.334    (0.001)&   0.333    (0.001)&  0.132    (0.041)\\
 &&b&      0.336     (0.002)&    0.333    (0.001)&    0.334     (0.001)&    0.334    (0.001)&   0.435    (0.010)\\
  \\
 & 3&
 a&      0.312     (0.014)&    0.334    (0.004)&    0.334    (0.004)&    0.333     (0.001)&  0.083      (0.063)\\
 &&b&      0.344    (0.003)&   0.333   (0.001)&    0.333    (0.001)&    0.333     (0.001)&   0.458     (0.016)\\
\hline\\

\end{tabular}
\caption{\small{Estimation results for Linear Case B, Gaussian error. Mean estimated values of the five estimators $\widehat{\theta}_{arma}$, $\widehat{\theta}_{N}$, $\widehat{\theta}_{SC}$, $\widehat{\theta}_{X}$ and $\widehat{\theta}_{naive}$ are presented for various values of $n$ ($1000$, 5000 or 10000) and s2n (0.5, 1.5, 3). True values are $a^0
  = 1/3$, $b^0 = 1/3$. MSEs are given in brackets.
}}
\label{table4}
\end{table}
}

\begin{figure}
\includegraphics
[width =10cm]{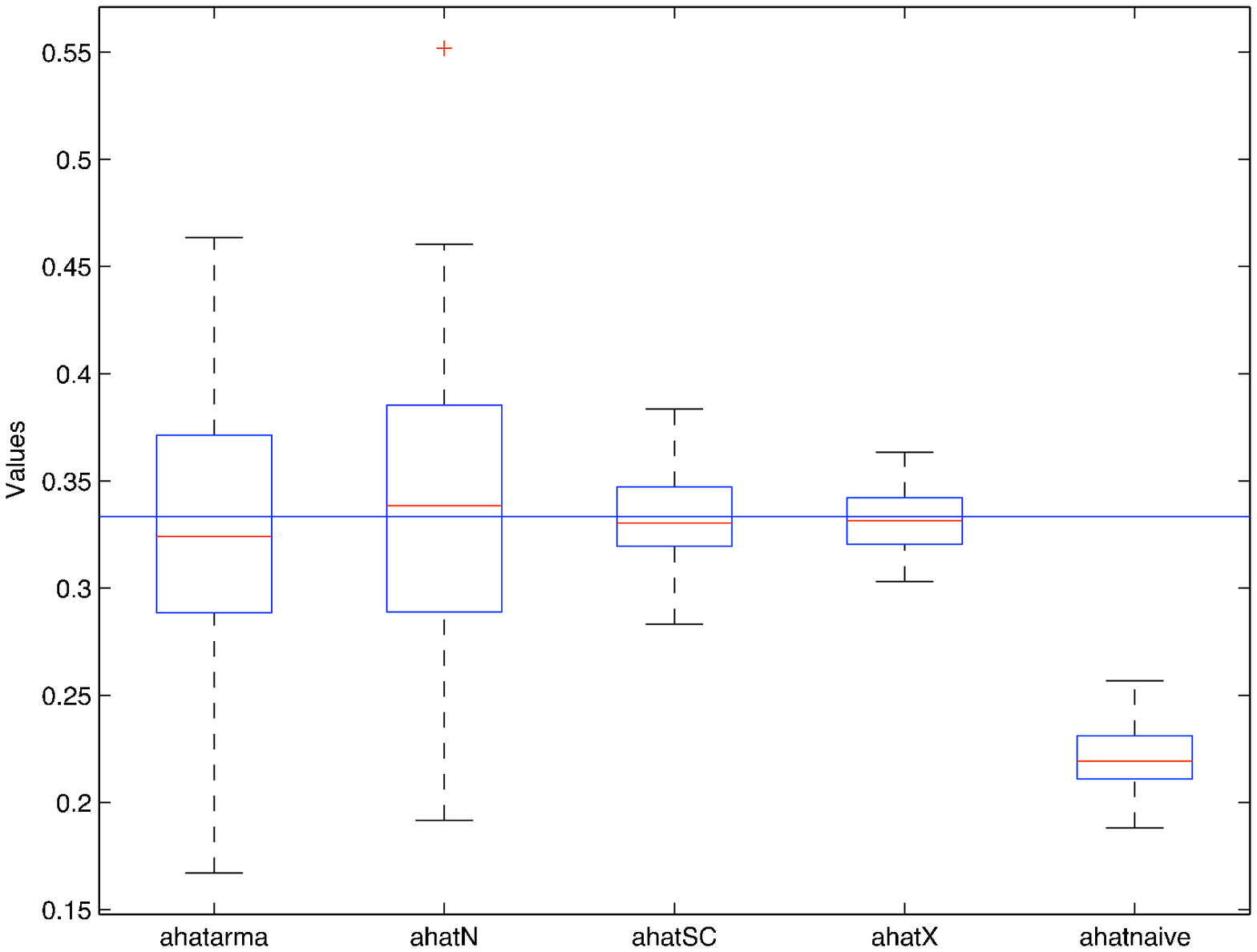}
\caption{\small{Results for linear Case B and  Gaussian error, with $n=5000$ and
  $\sigma_\varepsilon^2/\mbox{Var}(X)=0.5$. Box plots of the five estimators $\widehat{a}_{arma}$, $\widehat{a}_{N}$, $\widehat{a}_{SC}$, $\widehat{a}_{X}$ and $\widehat{a}_{naive}$, from left to right,  based on 100 replications. True value is $1/3$ (horizontal line). }}\label{fig0}
\end{figure}
\begin{figure}
\includegraphics
[width =10cm]{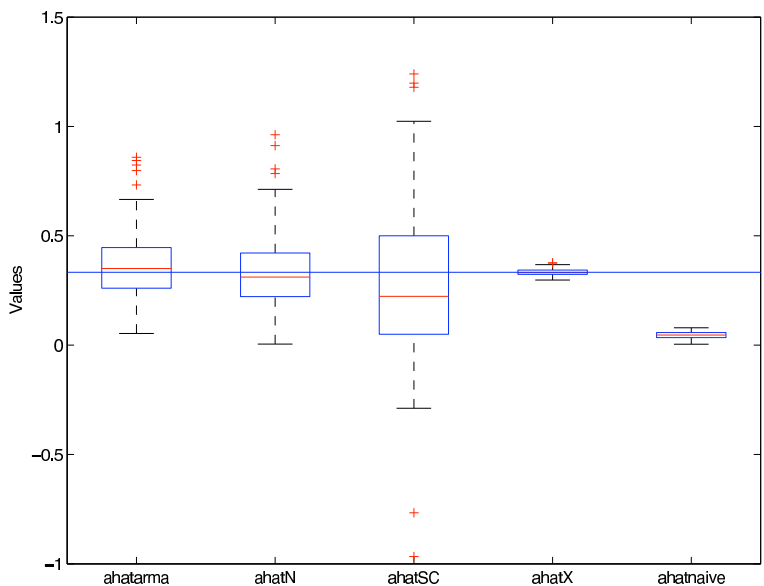}
\caption{\small{Results for linear Case B and  Gaussian error, with $n=5000$ and
  $\sigma_\varepsilon^2/\mbox{Var}(X)=6$. Box plots of the five estimators $\widehat{a}_{arma}$, $\widehat{a}_{N}$, $\widehat{a}_{SC}$, $\widehat{a}_{X}$ and $\widehat{a}_{naive}$, from left to right,  based on 100 replications. True value is $1/3$ (horizontal line).  }}\label{fig00}
\end{figure}


The first thing to notice is that, not surprisingly, $\widehat{\theta}_{naive}$ presents a
bias, whatever the values of $n$, $s2n$ and the error distribution.
The estimator $\widehat{\theta}_{X}$  has the good expected properties
(unbiased and  small MSE), but it  is based on the  observation of the
$X_i$'s.
The  previously  known  estimator $\widehat{\theta}_{arma}$  has  good
asymptotic properties. However its bias is often larger than the biases of
$\widehat{\theta}_{N}$  and $\widehat{\theta}_{SC}$,  except when $s2n=0.5$
and $\varepsilon$ is Gaussian.

We now consider the two estimators $\widehat{\theta}_{N}$
and $\widehat{\theta}_{SC}$.  Recall that their  construction requires
the choice of $w$.  Note first that, whatever the weight function $w$, the two estimators
$\widehat{\theta}_{N}$
and $\widehat{\theta}_{SC}$ present good convergence properties. Their biases and
MSEs decrease when $n$ increases.  When compared one to another, we can see that their
numerical behaviors are not the same. Namely for not too large $s2n$, $\widehat{\theta}_{SC}$ has a MSE smaller than
$\widehat{\theta}_{N}$    (see   Figure  \ref{fig0}   and  Tables
\ref{table1}-\ref{table4}, when $s2n\leq 3$). With large $s2n$,
the estimator  $\widehat{\theta}_{N}$ seems to  have better properties
(see Figure \ref{fig00} when $s2n=6$). This is expected since $N$ depends
on $\sigma_\varepsilon^2$  and is thus more sensitive  to small values
of $\sigma_\varepsilon^2$.
The error  distribution seems to have a slight infuence on the MSEs of
the two estimators.
The MSEs are often smaller when $f_\varepsilon$ is
the  Laplace density.  This may  be related  with the  theoretical
properties in density deconvolution.  In that context it is well
known that the  rate of convergence is slower  when $f_\varepsilon$ is
the Gaussian density.
The two  estimators $\widehat{\theta}_{N}$ and $\widehat{\theta}_{SC}$
have comparable  numerical behaviors  in the two  linear autoregressive
models. Let us  recall that in both cases, the  simulated chain $X$ are
non-mixing but are $\tau$-dependent. In Case A, the stationary distribution of
$X$ is continuous whereas it is not the case in Case B.
This explains the relative bad properties of $\widehat{\theta}_{arma}$
in Case B. Indeed, due to its construction, this estimator is expected to have good properties
when the stationary distribution of the Markov Chain is close to the Gaussian
distribution.  On the contrary our estimators
have similar behavior in both cases.

\subsection{Cauchy regression model}
\label{cauchy1}
We   consider   the   model  \eref{model}   with   $f_\theta(x)=\theta
/(1+x^2)=\theta f(x)$. The true parameter is $\theta^0=1.5$.
For the law of $\xi_0$  we take
$\xi_0\sim\mathcal{N}(0,0.01)$.  In this case, an empirical study shows that
$\sigma_X^2$ is about 0.1. Moreover  $\alpha_{{\bf X}}(k)=O(\kappa^k)$
for some $\kappa  \in ]0, 1[$ and the  Markov chain is $\alpha$-mixing
(see Appendix \ref{mixprop}).  For $w$ suitably chosen,
Theorem \ref{NAalpha} applies and states that $\hat{\theta}$ is asymptotically normal.
For the simulation, we start with $X_0$ uniformly distributed over $[0, 1]$, and we consider that the chain is close
to  the   stationary  chain  after   1000  iterations.  We   then  set
$X_i=X_{i+1000}$.

To our knowledge, the estimator $\widehat{\theta}$
is  the first consistent estimator in the literature for this regression function.
We first detail the estimator for two choices of the weight function $w$. Then we recall the classic estimator when $X$ is directly observed  and the so-called naive estimator.

\subsubsection{Expression of the estimator}
We consider the two following weight functions:
\begin{eqnarray}
\label{weightcauchy}
N_c(x)=(1+x^2)^2 \exp\{-x^2/(4 \sigma_\varepsilon^2)\} \mbox{ and }
SC_c(x) =(1+x^2)^2 \frac{1}{2\pi}\Big(\frac{2*\sin(x)}{x}\Big)^4,\end{eqnarray}
with
$\sigma_\varepsilon^2$ the variance of $\varepsilon$.
This  choice of  $w$ ensures  that  Conditions \eref{C11}-\eref{C11d4}
hold  and  our  method  allows  to  achieve  the  parametric  rate  of
convergence. As in the linear  case, these two weight functions differ
by their dependence on $\sigma_\varepsilon^2$ and their smoothness properties.
The  two associated  estimators  are  based  on the calculation of   $S_n(\theta)$, which can be
written as
$$S_n(\theta)=\frac{1}{n}\sum_{k=1}^n[Z_k^2I_w(Z_{k-1})+\theta^2I_{wf^2}(Z_{k-1})-2\theta
Z_k I_{wf}(Z_{k-1})],$$
where
\begin{eqnarray*}
I_{w}(Z)&=&\frac{1}{2\pi}\mathbb{R}e\int
(                  w)^*(u)\frac{e^{-iuZ}}{f_\varepsilon^*(-u)}du,\qquad
I_{wf}(Z)=\frac{1}{2\pi}\mathbb{R}e\int (w f)^*(u)\frac{e^{-iuZ}}{f_\varepsilon^*(-u)}du\\
\mbox{and} \quad
I_{wf^2}(Z)&=&\frac{1}{2\pi}\mathbb{R}e\int (w f^2 )^*(u)\frac{e^{-iuZ}}{f_\varepsilon^*(-u)}du.
\end{eqnarray*}
The estimator can be expressed as
\begin{eqnarray}
\label{coefcauchy}
\widehat{\theta}=\frac{\sum_{k=1}^n    Z_kI_{wf}(Z_{k-1})}{\sum_{k=1}^n
  I_{wf^2}(Z_{k-1})}.
\end{eqnarray}
In the following we  denote by $I_{wf,N_c}(Z)$,
$I_{wf^2,N_c}(Z)$, $I_{wf,SC_c}(Z)$ and $I_{wf^2,SC_c}(Z)$ respectively, the previous
integrals when the  weight function is either $w=N_c$  or $w=SC_c$. In the
same     way    we     denote     by    $\widehat{\theta}_{N_c}$     and
$\widehat{\theta}_{SC_c}$ the corresponding estimators of $\theta^0$.\\

\noindent $\bullet$ When $w=N_c$,
Fourier calculations provide that
\begin{eqnarray*}
 (N_cf)^*(t)&=&\sqrt{2\pi}\sqrt{2\sigma_\varepsilon^2}\exp(-\sigma_\varepsilon^2 t^2)\big(1+2\sigma_\varepsilon^2(1-2\sigma_\varepsilon^2t^2)  \big)\\
\mbox{ and }
 (N_cf^2)^*(t)&=&\sqrt{2\pi}\sqrt{2\sigma_\varepsilon^2}\exp(-\sigma_\varepsilon^2 t^2).
\end{eqnarray*}
Now, we can calculate
the integrals $I_{wf,N_c}(Z)$ and
$I_{wf^2,N_c}(Z)$.

\noindent 
If $f_\varepsilon$ is the Laplace distribution \eref{doublexp}, replacing
$f_\varepsilon^*$ by its expression we obtain
\begin{eqnarray*}
I_{wf,N_c}(Z)&=&\exp(-Z^2/(4\sigma^2_{\varepsilon}))\left[
Z^4-18Z^2\sigma_{\varepsilon}^2+Z^2+8\sigma_{\varepsilon}^4-10\sigma_{\varepsilon}^2\right]/(8\sigma_{\varepsilon}^2),\\
\mbox{ and }
I_{wf^2,N_c}(Z)&=&\exp(-Z^2/(4\sigma^2_{\varepsilon}))\big[1+\frac{1}{4}\big(1-\frac{Z^2}{2 \sigma^2_{\varepsilon}}\big)\big].
\end{eqnarray*}

\noindent 
If $f_\varepsilon$ is the Gaussian distribution \eref{gauss}, replacing
$f_\varepsilon^*$ by its expression we obtain
$$I_{wf,N_c}(Z)=\sqrt{2}e^{-Z^2/(2\sigma_\varepsilon^2)}(1-2\sigma_\varepsilon^2+4Z^2), \ \text{and} \quad I_{wf^2,N_c}(Z)=\sqrt{2}e^{-Z^2/(2\sigma_\varepsilon^2)}.$$

\noindent $\bullet$ When $w=SC_c$, easy calculations show that
$$I_{wf,SC_c}(Z)=I_{0,SC}(Z)+I_{2,SC}(Z)        \mbox{       and       }
I_{wf^2,SC_c}(Z)=I_{0,SC}(Z),$$
where $I_{0,SC}(Z)$ and $I_{2,SC}(Z)$ are defined by \eref{ISClin}. As explained before, the integrals $I_{0,SC}(Z)$ and $I_{2,SC}(Z)$ have no explicit form, whatever the error distributions, and are numerically approximated via the IFFT function.

\medskip
\subsubsection{Comparison with classical estimators.} We compare our estimators with two classical estimators, the usual least square estimator without observation noise, and the naive estimator.

\noindent $\bullet$ {Estimator  without noise.}
When $\varepsilon_i=0$, that is $(X_0, \ldots, X_n)$ is observed without errors, the parameter can be easily estimated
by the usual least square estimator
$$\widehat{\theta}_{X}=\frac{\sum_{i=1}^n X_if(X_{i-1})  }{\sum_{i=1}^n f^2(X_{i-1})}.$$

\medskip

\noindent $\bullet$  {Naive estimator.}
The idea  for the construction of  the naive estimator  is to replace
the unobserved $X_i$ by the observation $Z_i$
in the expression of $\widehat{\theta}_{X}$ to get
$$\widehat{\theta}_{naive}=\frac{\sum_{i=1}^n Z_if(Z_{i-1})  }{\sum_{i=1}^n f^2(Z_{i-1})}.$$
Classical results show that $\widehat{\theta}_{naive}$ is
an asymptotically biased estimator of $\theta^0$, which is confirmed by the simulation study.

\subsubsection{Simulations results}
\label{simu_results}
For each  error distribution, we  simulate 100   samples with  size  $n$, $n=500$,
$5000$ and $10000$. We consider different values of  $\sigma_\varepsilon$  such that the ratio signal
to noise $s2n=\sigma_\varepsilon^2/\mbox{Var}(X)$
is $0.5,1.5$ or $3$.

The comparison of the four estimators
is  based on the bias, the Mean Squared
Error (MSE), and the box plots.
The results are presented in Figure \ref{fig12} and Tables \ref{table5}-\ref{table6}.

The first thing to notice is that, not surprisingly, $\widehat{\theta}_{naive}$ presents a
bias,   whatever   the  values   of   $n$,   $s2n$   and  the   errors
distribution.  Moreover  it  converges   to  (false)  values  which are
different  according to  $s2n$ (see  Tables \eref{table5}-\eref{table6}).

The estimator $\widehat{\theta}_{X}$  has the good expected properties
(unbiased and  small MSE), but it  is based on the  observation of the
$X_i$'s.

We now compare  our two estimators illustrating the  influence of $w$,
$s2n$ and $f_\varepsilon$.
Globally, whatever the weight function $w$, the two estimators
$\widehat{\theta}$ present good convergence properties. Their biases and
MSEs decrease when $n$ increases.
The    MSEs   of    $\widehat{\theta}_{SC_c}$   increase    when   $s2n$
increases. This is not the case for the MSE of $\widehat{\theta}_{N_c}$. This is probably due to the fact that the
weight function  chosen for the  construction of $\widehat{\theta}_{N_c}$
depends  on  $\sigma_\varepsilon^2$.   This  estimator  is  thus  more
adaptive to changes in $s2n$.

{\tiny
\begin{table}
\begin{tabular}{llrrrr}
\hline
  & ratio&&\multicolumn{3}{c}{Estimator}\\
n&s2n&        \multicolumn{1}{c}{$\widehat{\theta}_{N_c}(MSE)$}         &
\multicolumn{1}{c}{$\widehat{\theta}_{SC_c}(MSE)$}&
\multicolumn{1}{c}{$\widehat{\theta}_X (MSE)$}
&\multicolumn{1}{c}{$\widehat{\theta}_{naive}(MSE)$}\\
\hline
$1000$ & 0.5& 1.5095 (0.0042)& 1.5024 (0.0006)& 1.5004 (0.0000)& 1.4333 (0.0050)\\
& 1.5& 1.5006 (0.0021)& 1.5005 (0.0013)& 1.5002 (0.0000)& 1.3657 (0.0190)\\
 &3& 1.5017 (0.0024) & 1.5005 (0.0024)& 1.5002 (0.0000)& 1.3267 (0.0314)
\\
\hline
$5000$& 0.5& 1.5045 (0.0008)& 1.5005 (0.0001)& 1.5003 (0.0000)& 1.4320
 (0.0047)\\
& 1.5& 1.5003 (0.0004)& 1.4994 (0.0003)& 1.4997 (0.0000)& 1.3647 (0.0185)\\
& 3& 1.4989 (0.0005)& 1.4992 (0.0005)& 1.5000 (0.0000)&
 1.3223 (0.0318)\\
\hline
$10000$& 0.5& 1.5033 (0.0004)& 1.5002 (0.0001)& 1.5000 (0.0000)&
 1.4315 (0.0047)\\
& 1.5& 1.5000 (0.0002)& 1.5000 (0.0001)& 1.4998 (0.0000)& 1.3650 (0.0183)\\
& 3& 1.4972 (0.0002)& 1.4970 (0.0002) & 1.4998 (0.0000)
 &1.3222 (0.0317)\\
\hline\\
\end{tabular}
\caption{\small{Estimation results for Cauchy, Laplace error. Mean estimated values of the four estimators $\widehat{\theta}_{N_c}$, $\widehat{\theta}_{SC_c}$, $\widehat{\theta}_{X}$ and $\widehat{\theta}_{naive}$ are presented for various values of  $n$ ($1000$,  5000 or 10000) and s2n (0.5, 1.5, 3). True value is $\theta^0=1.5$. MSE are given in brackets.}
}
\label{table5}
\end{table}}

{\tiny
\begin{table}
\begin{tabular}{llrrrr}
\hline
  & ratio&\multicolumn{3}{c}{Estimator}\\
n&s2n&         \multicolumn{1}{c}{$\widehat{\theta}_{N_c}(MSE)$}         &
\multicolumn{1}{c}{$\widehat{\theta}_{SC_c}(MSE)$}&
\multicolumn{1}{c}{$\widehat{\theta}_X (MSE)$}
&\multicolumn{1}{c}{$\widehat{\theta}_{naive}(MSE)$}\\
\hline
$1000$ & 0.5& 1.4979 (0.0027)& 1.4998 (0.0006)& 1.5000 (0.0000)& 1.4230 (0.0064)
\\
& 1.5& 1.4995 (0.0029)& 1.5001 (0.0015)& 1.5005 (0.0000)& 1.3336 (0.0287)\\
& 3& 1.5080 (0.0049)& 1.5058 (0.0042)& 1.4997 (0.0000)& 1.2832 (0.0487)
\\
\hline
$5000$  &  0.5& 1.5033  (0.0006)&  1.5011  (0.0001)& 1.4999  (0.0000)&
1.4250 (0.0057)\\
& 1.5& 1.5011  (0.0004) &  1.5001  (0.0003)& 1.4999  (0.0000)&
 1.3351 (0.0274)\\
 & 3 & 1.4998 (0.0009)& 1.4996 (0.0008)& 1.5002 (0.0000)& 1.2767
(0.0501)\\
\hline
$10000$ & 0.5&
    1.5017 (0.0003)& 1.4997 (0.0000)& 1.4996 (0.0000)& 1.4236 (0.0059)\\
& 1.5& 1.5025  (0.0003)& 1.5027  (0.0002)& 1.5001  (0.0000)&
1.3375 (0.0265)\\
 & 3& 1.5016 (0.0004)& 1.5021 (0.0004)& 1.5002 (0.0000)& 1.2778 (0.0495)
 \\
\hline\\

\end{tabular}
\caption{\small{Estimation results for Cauchy, Gaussian error. Mean estimated values of the four estimators $\widehat{\theta}_{N_c}$, $\widehat{\theta}_{SC_c}$, $\widehat{\theta}_{X}$ and $\widehat{\theta}_{naive}$ are presented for various values of  $n$ ($1000$,  5000 or 10000) and s2n (0.5, 1.5, 3). True value is $\theta^0=1.5$. MSE are given in brackets.}
}
\label{table6}
\end{table}}

\begin{figure}
\includegraphics
[width=10cm]{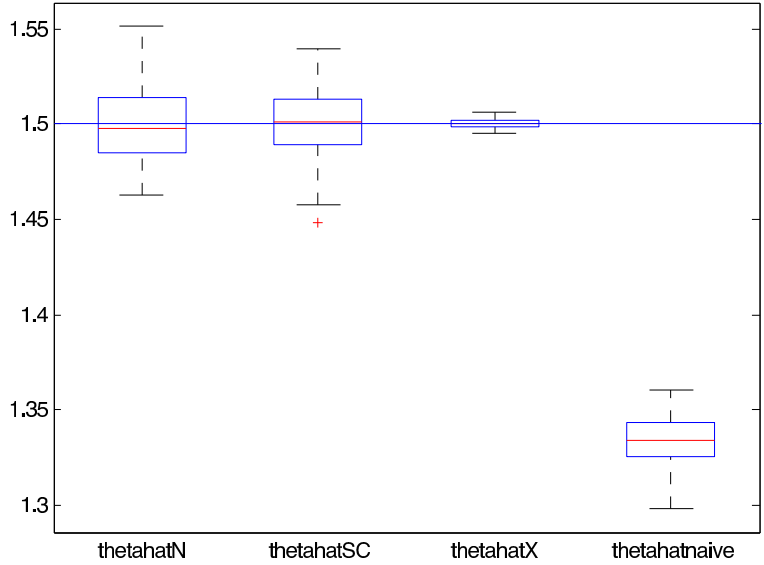}
\caption{\small{Results for Cauchy and  Gaussian error, with $n=5000$ and
  $\sigma_\varepsilon^2/\mbox{Var}(X)=1.5$. Box plots of the four estimators $\widehat{\theta}_{N_c}$, $\widehat{\theta}_{SC_c}$, $\widehat{\theta}_{X}$ and $\widehat{\theta}_{naive}$, from left to right,  based on 100 replications. True value is $1.5$ (horizontal line). }}\label{fig12}
\end{figure}

\section{A more general estimator}\label{generalcrit}
For a large number of regression functions,  a weight function $w$ such as the one involved in the definition of the estimator $\hat{\theta}$
can  be easily  exhibited.  Nevertheless  for some  specific regression
functions, it seems not straightforward to find a weight function such that  $(wf_\theta)^*/f_\varepsilon^*$ and
$(wf_\theta^2)^*/f_\varepsilon^*$ are integrable. We
refer to  Butucea and Taupin \citeyear{buttaupin} for  a more complete
discussion on this subject. Therefore, we propose a generalization of this estimator to relax these conditions.

\subsection{Definition of the general estimator}

The key idea for this construction is the following. We introduce
a density deconvolution kernel $K_{n,C_n}$  defined via its Fourier
transform $K_{n,C_n}^*$ by
\begin{eqnarray}
\label{kerdeconv}
K_{n,C_n}^*(t)=\frac{K^*(t/C_n)}{f_\varepsilon^*(-t)}:=\frac{K_{C_n}^*(t)}{f_\varepsilon^*(-t)},
\end{eqnarray}
where $K^*$  is the Fourier transform of  a kernel $K$ and  $C_n$ is a
sequence which tends to infinity with $n$. The kernel $K$  belongs to
 ${\mathbb L}^2({\mathbb R})$. Its Fourier transform $K^*$ is compactly supported
and satisfies $\vert 1-K^*(t)\vert\leq \ind_{\vert t\vert \geq 1}$.
Then, for any integrable
function
$\Phi$, one has  $\lim_{n\to\infty}n^{-1}\sum_{i=1}^n\Phi\star
K_{n,C_n}(Z_i)=\mathbb{E}(\Phi(X))$.
Hence we estimate $\mathbb{E}(\Phi(X))$ by $n^{-1}\sum_{i=1}^n\Phi\star
K_{n,C_n}(Z_i)$ instead of $n^{-1}\sum_{i=1}^n\Phi(X_i)$ which is not
available.
We then propose to estimate $ S_{\theta^0,P_X}(\theta)$ by
\begin{multline}
\label{Sngen}
 S_{n}(\theta)
=\frac{1}{n}\sum_{i=1}^n \mathbb{R}e\left[
\big((Z_i-f_\theta)^2w\big)\star K_{n,C_n}(Z_{i-1})\right]
\\=\frac{1}{n}\sum_{i=1}^n
\mathbb{R}e \!\!\int \!\!\left(Z_i-f_{\theta}(x)\right)^2w(x)\, K_{n,C_n}(Z_{i-1}-x)dx.
\end{multline}

Using this more general empirical criterion we propose to estimate $\theta^0$ by
\begin{eqnarray}
\widehat{\theta}=\arg\min_{\theta \in
\Theta}
S_{n}(\theta).\label{thetacgen}
\end{eqnarray}

Note that the general construction relies to a truncation of integrals
in  \eref{Sn}. Also note  that this  general construction  still works
under  Conditions   \eref{C11}-\eref{C11d4}.  It  suffices   to  chose
$K^*(t/C_n)=\ind_{\vert t \vert \leq C_n}$ with $C_n=+\infty$.

\subsection{Asymptotic properties under general assumptions}

This section presents the asymptotic properties of $\widehat {\theta}$
defined  by  \eref{thetacgen}  under  milder conditions  than  conditions
\eref{C11}-\eref{C11d4}, when one cannot exhibit a weight function $w$ ensuring that
these conditions  hold.  In  this  context  the  estimator  is  still
consistent, but  with a rate  which is not necessarily  the parametric
rate.
For  the  sake  of simplicity  we     only  consider the  case  of
$\alpha$-mixing Markov chains.

We assume that
\begin{align}
& \hypothese{emotilde2} \mbox{On }\Theta^0, \mbox{ the quantity }w^2(X_0)(Z_1-f_{\theta}(X_0))^{4}
\mbox{ and the absolute values of its derivatives} \\
&\notag \mbox{ with respect to } \theta \mbox{ up to order 2 } \mbox{ have a finite expectation.}
\\
&\hypothese{hypothese}\mbox{ The quantity } \sup_{n}\sup_{j\in \{1,\cdots, d\}}\mathbb{E}\Big(\sup_{\theta\in \Theta^\circ}\Big |
\frac{\partial}{\partial    \theta_j}   S_{n}(\theta)\Big    |   \Big)
\mbox{ is finite. }
\\
&\hypothese{cw1tilde} \sup_{\theta\in \Theta} \vert w f_\theta\vert,
\,\vert w\vert \mbox{ and } \sup_{\theta\in
  \Theta} \vert w
f^2_\theta\vert \mbox{ belong to }\mathbb{L}_1(\mathbb{R}).
\end{align}
We say that a function $\psi\in \mathbb{L}_1(\mathbb{R})$ satisfies \eref{condcons1}
if for a sequence $C_n$ we have
\begin{eqnarray}
\label{condcons1}
\min_{q=1,2}\parallel \psi^*(K_{C_n}^*-1)\parallel_q^2
+n^{-1}\min_{q=1,2}\left\| \frac{\psi^*K_{C_n}^*}{f_\varepsilon^*}\right\|_q^2
=o(1).
\end{eqnarray}

\begin{theo}
\label{thc1t}
Under the assumptions
\eref{Spttilde},     \eref{emo2tilde},    \eref{fepsnn},    \eref{df}
\eref{emotilde2} -
\eref{cw1tilde}, let $\widehat{\theta}$ be defined by \eref{thetacgen}
with $C_n$ such that \eref{condcons1} holds for $w$, $wf_\theta$
and $wf_\theta^2$ and their first derivatives with respect to $\theta$.
Assume that
the sequence $(X_k)$ is $\alpha$-mixing that is
$$\alpha_{\bf X}(k)\cv 0, \mbox{ as
} k \cv \infty.$$
Then
$\mathbb{E}(\|\widehat{\theta}-\theta^0\|_{\ell^2}^2)=o(1),$
 as $n\to \infty$ and $\widehat{\theta}$ is a consistent
 estimator of $\theta^0.$

\end{theo}
We now give upper bounds for the rates of convergence under two different types of assumptions:
\begin{align}
&\hypothese{R2} X_0\mbox{ admits a density } f_X \mbox{ with respect to the
  Lebesgue measure and there exist two } \\
& \nonumber \mbox{
  constants }
 C_1(f_{\theta^0}^2) \mbox{ and } C_2(f_{\theta^0})
\mbox{ such that }
\parallel f_{\theta^0}f_X\parallel_2^2\leq
C_1(f_{\theta^0}),\mbox{ and }\\&\nonumber \parallel f_{\theta^0}^2 f_X\parallel_2^2\leq
C_2(f_{\theta^0}^2).\\
&\hypothese{R3}\sup_{z\in
  \mathbb{R}}\mathbb{E}[f_{\theta^0}^2(X_0)f_\varepsilon(z-X_0)]
\mbox{ and } \sup_{z\in
  \mathbb{R}}\mathbb{E}[f_\varepsilon(z-X_0)]
\mbox{ are finite.}
\end{align}
These two  assumptions are mostly required for  technical reasons. The
following theorem still holds when  $X_0$ does not admit a density,
under a slightly different moment assumption.

\begin{theo}
\label{thv1t}
Suppose that  the assumptions of Theorem \ref{thc1t} hold. Assume moreover that the sequence
$(X_k)_{k\geq 0}$ is $\alpha$-mixing with
$\sum_{k\geq 1}\sqrt{\alpha_{\bf X}(k)}<\infty$,
and
that, for all $\theta\in \Theta$, the functions $w$, $f_\theta
w$ and $f_\theta^2w$ and their derivatives up to order 3 with respect to
$\theta$ satisfy \eref{condcons1}.

\textbf{1)} Assume that the sequence $X_0$ admits a
density  with respect  to  the Lebesgue  measure  and that  Assumption
\eref{R2}  holds.
Then
$\widehat{\theta}-\theta^0=O_p(\varphi_n^2)$ with
$\varphi_n=\|(\varphi_{n,j})\|_{\ell^2}$,
$\varphi_{n,j}^2= B_{n,j}^2+ V_{n,j}/n$, $j=1\ldots,d$, where
\begin{eqnarray*}
 B_{n,j}\!\!
&=&\!\!\min\left\lbrace  B_{n,j}^{[1]}, B_{n,j}^{[2]}\right\rbrace \mbox{ and
} V_{n,j}\!\!=\!\!\min\left\lbrace
  V_{n,j}^{[1]},
  V_{n,j}^{[2]}\right\rbrace
\end{eqnarray*}
and for $q=1,2$
\begin{eqnarray*}
 B_{n,j}^{[q]}&=&
\left\|(wf^{(1)}_{\theta,j})^*(K_{C_n}^*-1)\right\|^2_q
+\left\| (wf_{\theta^0}f^{(1)}_{\theta^0,j})^*(K_{C_n}^*-1)\right\|^2_q,
\end{eqnarray*}
and
\begin{eqnarray*}
 V_{n,j}^{[q]}\!\!&=&
\left\| (wf^{(1)}_{\theta^0,j})^*\frac{K_{C_n}^*}{f_\varepsilon^*}
\right\|_q^2+\left\|
(wf_{\theta^0}f^{(1)}_{\theta^0,j})^*\frac{K_{C_n}^*}{f_\varepsilon^*}
\right\|_q^2.
\end{eqnarray*}

\textbf{2)}
Assume that \eref{R3} holds. Then
$\widehat{\theta}-\theta^0=O_p(\varphi_n^2)$ with
$\varphi_n=\|(\varphi_{n,j})\|_{\ell^2}$,
$\varphi_{n,j}^2= B_{n,j}^2+ V_{n,j}/n$, $j=1\ldots,d$, where
$ B_{n,j}
=B_{n,j}^{[1]} \mbox{ and
} V_{n,j}=\min\left\lbrace
  V_{n,j}^{[1]},
  V_{n,j}^{[2]}\right\rbrace.$

\end{theo}

This theorem states an upper bound for the quadratic risk under very general
conditions.  It holds  under mild  conditions on  $w$,  $f_\theta$ and
$f_\varepsilon$.  We  refer to  Table  1 in  Butucea and  Taupin
\citeyear{buttaupin} for more details on the resulting rates.

\appendix

\section{Properties of the dependence coefficients and examples}

\subsection{Covariance inequalities and coupling}
\label{cov}
The following results are
the key arguments to prove the asymptotic normality of $\widehat{\theta}$.
We keep the same notations as in Definition \ref{defi}.

We first  recall a covariance inequality due to Rio \citeyear{Rio93}. For any positive random variable $Z$, let $Q_Z$ be the inverse cadlag of the tail function
$t
\to {\mathbb P}(Z>t)$.
Let $X$ and  $Y$ be two real valued random variables  such that
$\text{Cov}(X, Y)$ is well defined.
The following inequality holds
\begin{equation}\label{covalpha}
  |\text {Cov}(Y, X)|\leq 4 \int_0^{\alpha(\sigma(Y), \sigma(X))} Q_{|X|}(u) Q_{|Y|}(u) du \, .
\end{equation}

Next, we recall the coupling properties of
$\tau$ (see Dedecker and Prieur \citeyear{JDCPCoeff}): enlarging $\Omega$ if necessary, there exists $X^*$ distributed as $X$ and independent of ${\mathcal M}$ such that
\begin{equation}\label{coupling}
\tau({\mathcal M}, X)= {\mathbb E}(\|X-X^*\|_{\mathbb B}) \, .
\end{equation}

\subsection {Dependence properties of autoregressive models}
\label{mixprop}
We recall here  the mixing properties of the  autoregressive models
$$
X_i=f_{\theta^0}(X_{i-1})+\xi_i,
$$
that have been described in particular in the papers by
Mokkadem \citeyear{Mokkadem85} and Ango-Nz\'e \citeyear{AngoNze}. For instance, assume that
\begin{itemize}
\item the law of $\xi_0$ has a density $f_\xi$ such that  $f_\xi>c>0$ on a neighborhood of zero,
and there exists $S \geq 1$ such that ${\mathbb E}(|\xi_0|^S) < \infty$.
\item $f_{\theta^0}$ is continuous and there exist $R\geq 1$ and $\rho \in ]0, 1[$ such that: for any $|x|\geq R$,
$|f_{\theta^0}(x)| \leq \rho |x|$.
\end{itemize}
Then there exists a unique invariant probability measure, and  the stationary  Markov chain $(X_i)_{i \geq 0}$ satisfies
$\alpha_{{\bf X}}(k)=O(\kappa^k)$ for any $\kappa \in ]\rho, 1[$ and is $\alpha$-mixing.

\noindent Now if the second point is weakened to
\begin{itemize}
\item $f_{\theta^0}$ is continuous and there exist $R\geq 1$ and $\delta \in ]0, 1[$ such that: for any $|x|\geq R$,
$|f_{\theta^0}(x)| \leq |x|(1-|x|^{-\delta})$.
\end{itemize}
Then there exists a unique invariant probability measure, and the stationary Markov chain $(X_i)_{i \geq 0}$ satisfies $\alpha_{{\bf X}}(k)=O(k^{1-S/\delta})$ and is $\alpha$-mixing.

Now, if  we do not assume that  $\xi_0$ has a density,  then the chain
may not be  $\alpha$-mixing (and not even irreducible). However, under appropriate assumptions on $f_{\theta^0}$, it is still
possible to obtain upper bounds for the coefficient $\tau$. For instance assume that
\begin{itemize}
\item  there exists $S \geq 1$ such that ${\mathbb E}(|\xi_0|^S) < \infty$.
\item  $|f_{\theta^0}(x)-f_{\theta^0}(y)|\leq \rho |x-y|$ for some $\rho \in ]0, 1[$.
\end{itemize}
Then there exists a unique invariant probability measure, and the stationary Markov chain $(X_i)_{i \geq 0}$ satisfies $\tau_{{\bf X}, 2}(k)=O(\rho^k)$ and is $\tau$-dependent.
Now if the second point is weakened to
\begin{itemize}
\item there exist $ \delta$ in $[0,1[$ and   $C$ in $]0,1]$ such that $|f'(t)| \leq 1 - C(1+|t|)^{-\delta}$ almost everywhere.
\end{itemize}
Then there exists a unique invariant probability measure, and for $S>1 + \delta$ the stationary Markov chain $(X_i)_{i \geq 0}$ satisfies $\tau_{{\bf X}, 2}(n)=O(n^{(\delta+1-S)/\delta})$ and is $\tau$-dependent.

\section{proofs of Theorems}
\label{proofsth}
\setcounter{equation}{0}

\subsection{Proof of Theorem \ref{consist}}

The main point of the proof consists in showing the two following points

i) for any $\theta$ in $\Theta$,
$S_{n}(\theta)\cvL1 S_{\theta^0,P_X}(\theta)$, with
$S_{\theta^0,P_X}(\theta)$ admitting a unique minimum in
$\theta=\theta^0$.

ii) For $\omega_2(n,\rho)$ defined as
$\omega_2(n,\rho)=\sup\left\lbrace \vert S_{n} (\theta)-S_{n}(\theta')\vert:\|
\theta-\theta'\|_{\ell^2}\leq \rho \right\rbrace,
$
there exists a sequence $\rho_k$ tending to 0, such
that
\begin{eqnarray}
{\mathbb E} (\omega_2(n,\rho_k))=O(\rho_k).\label{equic2}
\end{eqnarray}
Let us start with the proof of i) by writing that
$$ S_{n}(\theta)=\frac{1}{n}\sum_{k=1}^n \Psi (Z_k,Z_{k-1}), \mbox{ with }\Psi(Z_1,Z_0)=\frac{1}{2\pi }
\mathbb{R}e\int \frac{\left(\big(Z_1-f_{\theta}\big)^2w\right)^*(t)e^{-itZ_0}}{f_\varepsilon^*(-t)}dt,$$
that is seen as a function of a strictly stationary and ergodic sequence
of   random  variables.   By  the   ergodic  theorem   and  Assumption
\eref{emotilde} we conclude that for any $\theta\in \Theta$,
$$S_{n}(\theta)\cvL1 \mathbb{E}(\psi(Z_1,Z_0))=S_{\theta^0,P_X}(\theta).$$

It remains now to check that there exists a sequence $\rho_k$ tending to 0, such that \eref{equic2} holds. This follows by the assumption \eref{C11d1} and by writing that
\begin{equation}\label{majoration}
\sup_{\parallel \theta-\theta^\prime\parallel_{\ell^2}\leq \rho }|S_{n}(\theta)-S_{n}(\theta^{\prime})|
\leq \sup_{\parallel \theta-\theta^\prime\parallel_{\ell^2}\leq \rho }\parallel \theta-\theta^\prime\parallel_{\ell^2}\sup_{\theta\in \Theta^0}\parallel  S^{(1)}_{n}(\theta)\parallel_{\ell^2}.
\end{equation}
\qed

\subsection{Proof of Theorem \ref{NAalpha}}

By using a Taylor expansion based on the smoothness
properties of $\theta\mapsto w f_\theta$
and the consistency of $\widehat{\theta}$, we obtain
$$
0=S^{(1)}_{n}(\widehat{\theta})=S^{(1)}_{n}(\theta^0)+
S^{(2)}_{n}(\theta^0)(\widehat{\theta}-\theta^0)+
R_{n}(\widehat{\theta}-\theta^0),$$
 with $R_{n}$ defined by
\begin{equation}
R_{n}=\int_{0}^{1}[S^{(2)}_{n}(\theta^0+s(\widehat{\theta}-\theta^0))
-S^{(2)}_{n}(\theta^0)]ds. \label{Rn}
\end{equation}
This implies that
\begin{equation}
\label{base} \widehat{\theta}-\theta^0=-[S^{(2)}_{n}(\theta^0)+R_{n}]^{-1}
S^{(1)}_{n}(\theta^0).
\end{equation}
Consequently, we have to check the three following points.
\begin{itemize}
\item[i)] $\sqrt{n}S^{(1)}_{n}(\theta^0)\cvl \mathcal{N}(0,\Sigma_{0,1})$;
\item[ii)] $ S^{(2)}_{n}(\theta^0) \cvp S^{(2)}_{\theta^0,P_X}(\theta^0);$
\item[iii)] $R_{n}$ defined in \eref{Rn} satisfies
  $R_{n}\cvp 0.$
\end{itemize}

Note that the covariance matrix $\Sigma_{0,1}$ in i) satisfies $\Sigma_{0,1}=\Sigma / 4 \pi^2$, with $\Sigma$ defined by the equation \eref{somcov} below.
Consequently, according to ii) and iii),
the covariance matrix  $\Sigma_{1}$  satisfies
\begin{equation}\label{secondequality}
\Sigma_{1}=\frac{1}{4 \pi^2}(S^{(2)}_{\theta^0,P_X}(\theta^0))^{-1}\Sigma (S^{(2)}_{\theta^0,P_X}(\theta^0))^{-1}, \quad \text{with $\Sigma$ defined by \eref{somcov}}.
\end{equation}

\medskip

\begin{center}
\textit{\textbf{Proof of i)}}
\end{center}
Under Assumption \eref{C11d1},
\begin{eqnarray*}
\left(\sqrt{n} S^{(1)}_{n}(\theta^0)\right)_i
=\frac{1}{2\pi\sqrt{n}}\sum_{k=1}^n \mathbb{R}e\int 
\left(\frac{\partial}{\partial \theta_i} ((Z_k-f_{\theta})^2)w \Big \vert_{\theta=\theta^0}\right)^*(t) \frac{e^{-itZ{k-1}}}{f_\varepsilon^*(-t)}dt.
  \end{eqnarray*}
  We have thus to prove that
  \begin{eqnarray*}
  \frac{1}{2\pi\sqrt{n}}\sum_{k=1}^n \mathbb{R}e\int \big( -2(Z_k-f_{\theta^0})f^{(1)}_{\theta^0}w\big)^*(t)
  \frac{e^{-itZ{k-1}}}{f_\varepsilon^*(-t)}dt\cvl \mathcal{N}(0,\Sigma_{0,1}).
  \end{eqnarray*}
We first use that $\mathbb{E}( S_{n}(\theta))=S_{\theta^0,P_X}(\theta)$ and thus
$\mathbb{E}(S^{(1)}_{n}(\theta^0))=S^{(1)}_{\theta^0,P_X}(\theta^0)=0.$
Next we write
$$
\sqrt{n}S^{(1)}_{n}(\theta^0)=\sqrt{n}S^{(1)}_{n}(\theta^0)-
\mathbb{E}[\sqrt{n} S^{(1)}_{n}(\theta^0)]=
\frac{1}{2\pi\sqrt{n}}\sum_{k=1}^n T_k
$$
with $T_k=-2W_{k,1}+2W_{k,2}$, and 
\begin{eqnarray*}
W_{k,1}&=&Z_k \mathbb{R}e\int \big(f^{(1)}_{\theta^0}w\big)^*(t)
  \frac{e^{-itZ_{k-1}}}{f_\varepsilon^*(-t)}dt- \mathbb{E}\left[ Z_k \mathbb{R}e\int \big(f^{(1)}_{\theta^0}w\big)^*(t)
  \frac{e^{-itZ_{k-1}}}{f_\varepsilon^*(-t)}dt  \right]\\
  W_{k,2}&=&\mathbb{R}e\int \big(f_{\theta^0}f^{(1)}_{\theta^0}w\big)^*(t)
  \frac{e^{-itZ_{k-1}}}{f_\varepsilon^*(-t)}dt-\mathbb{E}\left[ \mathbb{R}e\int \big(f_{\theta^0} f^{(1)}_{\theta^0}
w\big)^*(t)
  \frac{e^{-itZ_{k-1}}}{f_\varepsilon^*(-t)}dt  \right].
\end{eqnarray*}
Let ${\mathcal M}_1=\sigma(X_0, X_1, \varepsilon_0, \varepsilon_1)$.
According to Dedecker and Rio \citeyear{DedeckerRio},
$n^{-1/2}\sum_{k=1}^n T_k$  converges to a centered Gaussian vector with
covariance matrix
\begin{eqnarray}\label{somcov}
\Sigma=\mathrm {Cov}(T_1, T_1)  + 2 \sum_{k>1} \mathrm {Cov}(T_1, T_k) \, ,
\end{eqnarray}
as soon as for any $ (p, q)$ in $\{1,\cdots,d\} \times \{1,\cdots,d\}$
\begin{eqnarray}
\label{TLC}
\sum_{k=3}^\infty\mathbb{E}\vert (T_1)_p\mathbb{E}((T_k)_q|\mathcal{M}_1)\vert <\infty.
\end{eqnarray}
For any $(p, q)$ in $\{1,\cdots,d\} \times \{1,\cdots,d\}$
and any $i,j  \in \{1,2\}$, we shall give an upper bound for
\begin{eqnarray*}
\mathbb{E}\left\vert (W_{1,i})_p \mathbb{E}((W_{k,j})_q|\mathcal{M}_1)\right\vert.
\end{eqnarray*}
We first notice that the sequence $(\varepsilon_k, \varepsilon_{k-1})$ is independent of $\mathcal{M}_1\vee \sigma(X_k, X_{k-1})$. It follows that for $i,j\in \{1,2\}$,
\begin{eqnarray*}
\mathbb{E}\left\vert (W_{1,i})_p \mathbb{E}((W_{k,j})_q|\mathcal{M}_1)\right\vert=
\mathbb{E}\left\vert (W_{1,i})_p \mathbb{E}((\tilde W_{k,j})_q|\mathcal{M}_1)\right\vert,
\end{eqnarray*}
with
\begin{eqnarray*}
(\tilde W_{k,1})_q&=&X_k \int \big( f^{(1)}_{\theta^0,q}w\big)^*(t)
  e^{-itX_{k-1}}dt- \mathbb{E}\left[ X_k \int \big(f^{(1)}_{\theta^0,q}w\big)^*(t)
  e^{-itX_{k-1}}dt  \right]\\
  (\tilde W_{k,2})_q&=&\int \big(f_{\theta^0} f^{(1)}_{\theta^0,q}w\big)^*(t)
  e^{-itX_{k-1}}dt-\mathbb{E}\left[\int \big(f_{\theta^0}  f^{(1)}_{\theta^0,q}w\big)^*(t)
  e^{-itX_{k-1}}dt  \right].
 \end{eqnarray*}
Next, since ${\mathbb P}_{(X_{k-1}, X_k)|\sigma(\varepsilon_0, \varepsilon_1, X_0, X_1)}={\mathbb P}_{(X_{k-1}, X_k)|\sigma(X_1)}$, we infer that
\begin{eqnarray*}
\mathbb{E}\left\vert (W_{1,i})_p \mathbb{E}((W_{k,j})_q|\mathcal{M}_1)\right\vert=
\mathbb{E}\left\vert (W_{1,i})_p \mathbb{E}((\tilde W_{k,j})_q|X_1)\right\vert.
\end{eqnarray*}
Next we use that under Condition \eref{C11d1},
\begin{eqnarray*}
\vert (W_{1,1})_p\vert &\leq& \vert Z_1\vert \int \left\vert \big(
 f^{(1)}_{\theta^0,p}w\big)^*(t)
  \frac{e^{-itZ_{0}}}{f_\varepsilon^*(-t)}\right\vert dt+
  \mathbb{E}\left\lbrace\vert Z_1\vert \int \left\vert \big(
 f^{(1)}_{\theta^0,p}w\big)^*(t)
  \frac{e^{-itZ_{0}}}{f_\varepsilon^*(-t)}\right\vert dt \right\rbrace\\
  &\leq &\vert Z_1\vert \int \left\vert \big(f^{(1)}_{\theta^0,p}w\big)^*(t)
  \frac{1}{f_\varepsilon^*(-t)}\right\vert dt+
  \mathbb{E}\left\lbrace\vert Z_1\vert \int \left\vert \big(
 f^{(1)}_{\theta^0,p}w\big)^*(t)
  \frac{1}{f_\varepsilon^*(-t)}\right\vert dt \right\rbrace\\
  &\leq& C_1(\vert Z_1\vert+\mathbb{E}(\vert Z_1\vert)).
  \end{eqnarray*}
  In the same way we get that
  $\vert (W_{1,2})_p\vert\leq C_2$.

Now, since $\varepsilon_1$ is independent of $X_1$, for $j \in \{1, 2\}$
\begin{eqnarray}\label{control1}\mathbb{E}\left\vert (W_{1,1})_p \mathbb{E}((\tilde W_{k,j})_q|X_1)\right\vert&\leq&
C_1\mathbb{E}\left[(\vert Z_1\vert+ {\mathbb E}(|Z_1|))\left\vert\mathbb{E}((\tilde W_{k,j})_q|X_1)\right\vert\right]\nonumber\\
&\leq & C\mathbb{E}\left[(\vert X_1\vert+ {\mathbb E}(|X_1|))\left\vert\mathbb{E}((\tilde W_{k,j})_q|X_1)\right\vert\right].
\end{eqnarray}
In the same way
\begin{equation}\label{control2}\mathbb{E}\left\vert (W_{1,2})_p \mathbb{E}((\tilde W_{k,j})_q|X_1)\right\vert\leq
C\mathbb{E}\left\vert\mathbb{E}((\tilde W_{k,j})_q|X_1)\right\vert.
\end{equation}
Note that
$$
\mathbb{E}\left[(\vert X_1\vert+ {\mathbb E}(|X_1|))\left\vert\mathbb{E}((\tilde W_{k,1})_q|X_1)\right\vert\right]=
\mathrm{Cov}((\vert X_1\vert+ {\mathbb E}(|X_1|))\text{sign}(\mathbb{E}((\tilde W_{k,1})_q|X_1)), ( \tilde W_{k, 1})_q).
$$
Now, we use the covariance inequality (\ref{covalpha}).
Note first that
$$(\vert X_1\vert+ {\mathbb E}(|X_1|))\text{sign}(\mathbb{E}((\tilde W_{k,1})_q|X_1))
\leq \vert X_1\vert+ {\mathbb E}(|X_1|)$$ and
$$
\vert (\tilde W_{1,1})_q\vert \leq D(\vert X_1\vert+ {\mathbb E}(|X_1|))\, .
$$
Since $(X_i)_{i \geq 0}$ is a strictly stationary Markov chain, it is well
known that
\begin{equation}\label{Mchains}
\alpha(\sigma(X_1), \sigma(X_{k-1}, X_k))=\alpha(\sigma(X_1), \sigma(X_{k-1}))=\alpha_{\bf X}(k-2)\, .
\end{equation}
 Hence, applying (\ref{covalpha}),
$$
\mathbb{E}\left\vert (W_{1,1})_p \mathbb{E}((\tilde W_{k,1})_q|X_1)\right\vert
\leq C
\int_0^{\alpha_{\bf X}(k-2)} Q^2_{|X_1|}(u)du \, .
$$
We conclude  that
\begin{eqnarray*}
\sum_{k\geq 3}\mathbb{E}\left\vert (W_{1,1})_p \mathbb{E}((W_{k,1})_q|\mathcal{M}_1)\right\vert\leq
C \sum_{k\geq 3}
\int_0^{\alpha_{\bf X}(k-2)}Q^2_{|X_1|}(u) du.
\end{eqnarray*}
Finally, using similar arguments for the three quantities $\sum_{k\geq
  3}\mathbb{E}\left\vert                                    (W_{1,2})_p
  \mathbb{E}((W_{k,1})_q|\mathcal{M}_1)\right\vert$, \\
$\sum_{k\geq 3}\mathbb{E}\left\vert (W_{1,1})_p \mathbb{E}((W_{k,2})_q|\mathcal{M}_1)\right\vert$
and        $\sum_{k\geq       3}\mathbb{E}\left\vert       (W_{1,2})_p
  \mathbb{E}((W_{k,2})_q|\mathcal{M}_1)\right\vert$
we conclude that
 $$\sqrt{n} S^{(1)}_{n}(\theta^0)\cvl \mathcal{N}(0,\Sigma/(4 \pi^2))$$ as soon as
$$\sum_{k\geq 1}
\int_0^{\alpha_{\bf X}(k)}Q^2_{|X_1|}(u) du<\infty.$$

\qed

\medskip

\begin{center}
\textit{\textbf{Proof of ii)}}
\end{center}
Under Condition \eref{C11d2}, for $j,k=1,\cdots,d$,
\begin{eqnarray}
\label{Snd2}
\left( S^{(2)}_{n}(\theta)\right)_{j,k}
=\frac{1}{2\pi n}\sum_{\ell=1}^n \mathbb{R}e\int \left(-2Z_\ell \frac{\partial^2}{\partial \theta_j \partial \theta_k} (f_{\theta}
    w)
      +\frac{\partial^2}{\partial \theta_j \partial\theta_k} (f_{\theta}^2w)
\right)^*(t)
      \frac{e^{-itZ_{\ell-1}}}{f_\varepsilon^*(-t)}dt
      \end{eqnarray} and by applying the ergodic theorem we get that
      $$ S^{(2)}_{n}(\theta^0)\cvp S^{(2)}_{\theta^0,P_X}(\theta^0).$$ \qed

\medskip

\begin{center}
\textit{\textbf{Proof of iii)}}
\end{center}
Starting from \eref{Rn} and  \eref{Snd2}, the point iii) follows from the assumption  \eref{C11d3} on the properties of the derivatives at order 3 of
$w f_\theta$ and $w f_\theta^2$. \qed

\subsection{Proof of Theorem \ref{NAtau}}
We follow the proof of Theorem \ref{NAalpha} and keep the same notations. We have to check that the condition (\ref{TLC}) holds.
We start from the inequalities (\ref{control1}) and (\ref{control2}). For clarity, let us write
$$
(\tilde W_{k,1})_q=(\tilde W_{k,1})_q(X_{k}, X_{k-1})\, .
$$
Let $\psi_M$ be the truncating function
defined by $\psi_M(x)=(x \wedge M)\vee (-M)$. Applying (\ref{coupling}), let $(X_k^*, X_{k-1}^*)$ be the random variable distributed as
$(X_k, X_{k-1})$ and independent of $X_1$ such that
$$
  \frac{1}{2}(\|X_k-X_k^*\|_1+\|X_{k-1}-X^*_{k-1}\|_1) = \tau(\sigma(X_1), (X_{k-1}, X_k))\leq \tau_{X,2}(k-2)\, .
$$
Define the constants $K_1$ and $K_2$ by
\begin{equation*}
K_1=\int \left |\big(f^{(1)}_{\theta^0,q}w\big)^*(t)\right| dt  < \infty \, , \quad
K_2=\int |t| \left |\big(f^{(1)}_{\theta^0,q}w\big)^*(t)\right| dt  < \infty \, .
\end{equation*}
Clearly
$$
|X_1 \mathbb{E}((\tilde W_{k,1})_q(X_k, X_{k-1})|X_1)|\leq M |\mathbb{E}((\tilde W_{k,1})_q(X_k, X_{k-1})|X_1)|+ K_2|X_1|{\mathbf 1}_{|X_1|>M}(|X_k|+{\mathbb E}(|X_k|))\, .
$$
Now, since $(X_k^*, X_{k-1}^*)$ is independent of $X_1$, one has that
$$
  |\mathbb{E}((\tilde W_{k,1})_q(X_k, X_{k-1})|X_1)|= |\mathbb{E}((\tilde W_{k,1})_q(X_k, X_{k-1})-(\tilde W_{k,1})_q(X^*_k, X^*_{k-1})|X_1)|\, .
$$
By definition of $(\tilde W_{k,1})_q(X_k, X_{k-1})$, there exists a constant $C$ such that
\begin{multline*}
|(\tilde W_{k,1})_q(X_k, X_{k-1})-(\tilde W_{k,1})_q(X^*_k, X^*_{k-1})-((\tilde W_{k,1})_q(\psi_M(X_k), X_{k-1})-(\tilde W_{k,1})_q(\psi_M(X^*_k), X^*_{k-1}))|\\\leq
 C(|X_k|{\mathbf 1}_{|X_k|>M}+ |X_k^*|{\mathbf 1}_{|X_k^*|>M}).
\end{multline*}
Hence
\begin{eqnarray*}
| \mathbb{E}((\tilde W_{k,1})_q(X_k, X_{k-1})|X_1)|& \leq &  | \mathbb{E}((\tilde W_{k,1})_q(\psi_M(X_k), X_{k-1})-(\tilde W_{k,1})_q(\psi_M(X^*_k), X^*_{k-1})|X_1)|\\
& + & C(|X_k|{\mathbf 1}_{|X_k|>M}+ |X_k^*|{\mathbf 1}_{|X_k^*|>M})\, .
\end{eqnarray*}
Since  $\psi_M$  is 1-Lipschitz  and  bounded  by  $M$, and  since  $x
\rightarrow  \exp(itx)$ is  $|t|$-Lipschitz  and bounded  by 1,  under
Condition \eref{C11d4},
one has
$$
|(\tilde W_{k,1})_q(\psi_M(X_k), X_{k-1})-(\tilde W_{k,1})_q(\psi_M(X^*_k), X^*_{k-1})| \leq MK_2 |X_{k-1}-X_{k-1}^*| + K_1 |X_k-X_k^*| \, .
$$
It follows that
\begin{eqnarray*}
|X_1 \mathbb{E}((\tilde W_{k,1})_q(X_k, X_{k-1})|X_1)| &\leq & K_2|X_1|{\mathbf 1}_{|X_1|>M}(|X_k|+{\mathbb E}(|X_k|))\\
&+& M^2K_2 |X_{k-1}-X_{k-1}^*| + MK_1 |X_k-X_k^*|)\\
&+& CM(|X_k|{\mathbf 1}_{|X_k|>M}+ |X_k^*|{\mathbf 1}_{|X_k^*|>M})\, .
\end{eqnarray*}
Using that
$$
|X_1|{\mathbf 1}_{|X_1|>M}|X_k| \leq \frac32 X_1^2{\mathbf 1}_{|X_1|>M} +\frac12 X_k^2{\mathbf 1}_{|X_k|>M} \, ,
$$
we infer from \eref{control1} with $j=1$ that there exists a positive constant $K$ such that
$$
\mathbb{E}\left[\left \vert (W_{1, 1})_p\mathbb{E}((\tilde W_{k,1})_q|X_1)\right\vert\right] \leq
K(L(M^2)+ M(M+1)\tau_{{\bf X},2}(k-2))\, ,
$$
where $L(t)={\mathbb E}(X_0^2{\mathbf 1}_{X_0^2>t})$.
Let then $G(t)=t^{-1}L(t)$, and let $G^{-1}$ be the inverse cadlag of $G$. Choose then $M^2=G^{-1}(\tau_{{\bf X},2}(k-2))$. We obtain that
\begin{multline*}
\mathbb{E}\left[\left\vert        (W_{1,       1})_p\mathbb{E}((\tilde
    W_{k,1})_q|X_1)\right\vert\right] \leq 2K(2G^{-1}(\tau_{{\bf X},2}(k-2))\tau_{{\bf X},2}(k-2)
+
\sqrt{G^{-1}(\tau_{{\bf X},2}(k-2))}\tau_{{\bf X},2}(k-2))\, .
\end{multline*}
It follows that
$$
  \sum_{k\geq 3} \mathbb{E}\left[\left\vert (W_{1, 1})_p\mathbb{E}((\tilde W_{k,1})_q|X_1)\right\vert\right] < \infty \quad \text{as soon as}
  \quad \sum_{k>0}G^{-1}(\tau_{{\bf X},2}(k)) \tau_{{\bf X},2}(k) < \infty \, .
$$
Easier control holds for the other  terms in (\ref{control1}) and (\ref{control2}). Consequently (\ref{TLC}) holds as soon as
(\ref{condtau}) holds,
and the proof is complete.

\subsection{Proof of Theorem \ref{thc1t}}

The  proof  of  the  consistency  under  the  assumptions  of  Theorem
\ref{thc1t} is quite different from the proof of the consistency under
Conditions \eref{C11}-\eref{C11d1} in Theorem \ref{consist}. This comes from the fact that $S_{n}(\theta)$ is now a triangular array of the form
$$S_{n}(\theta)=\frac{1}{n}\sum_{k=1}^n \Psi_n(Z_k,Z_{k-1}) \mbox{ with }\Psi_n(Z_1,Z_{0})=\frac{1}{2\pi}\mathbb{R}e\int \frac{\left(\big(Z_1-f_{\theta}\big)^2w\right)^*(t)e^{-itZ_0}K_{C_n}^*(t)}{f_\varepsilon^*(-t)}dt.
$$
In this context we show that

i) For all $\theta$ in $\Theta$,
$\mathbb{E}[(S_{n}(\theta)-S_{\theta^0,P_X}(\theta))^2]=o(1)$ as $n \to \infty.$

ii) The control \eref{equic2} holds.

Note first that ii) follows from the upper bound \eref{majoration} and Assumption \eref{hypothese}.

For the proof of i) we check that for all $\theta \in \Theta,$
\begin{equation}\label{cons}
\mathbb{E}[S_{n} (\theta)] - S_{\theta^0,P_X}(\theta) = o(1)
\quad  \mbox{ and }
\mbox{Var}(S_{n}(\theta))= o(1),
\quad \text{ as } n \to \infty.
\end{equation}

\proof[Proof of the first part of $\eref{cons}$]

Since $Z_0=X_0+ \varepsilon_0$, with $\varepsilon_0$ independent of $(Z_1, X_0)$, it follows that
\begin{eqnarray*}
\mathbb{E}[S_{n} (\theta)]
=\mathbb{E} \left[\mathbb{R}e\left( (Z_1 - f_\theta)^2
    \,w\right)
\star K_{n,C_n}(Z_0)\right]=\mathbb{E} \left[\left( (Z_1 - f_\theta)^2
    \,w\right)
\star K_{C_n}(X_0)\right]\, ,
\end{eqnarray*}
hence
\begin{eqnarray*}
\mathbb{E}[S_{n} (\theta)] - S_{\theta^0,P_X}(\theta)
&=&\frac{1}{2\pi}\iint (f^2_{\theta^0}(x)+\sigma^2_\xi + \sigma^2_\varepsilon) e^{-iux} w^*(u)(K_{C_n}^*-1)(u) du
P_X(dx)\\&&-\frac{1}{\pi}\iint f_{\theta^0}(x) e^{-iux}
(f_\theta w)^*(u)(K_{C_n}^*-1)(u)du P_X(dx)\\
&&+\frac{1}{2\pi}\iint e^{-iux}
(f^2_{\theta}w)^*(u)(K_{C_n}^*-1)(u)P_X(dx)du.
\end{eqnarray*}
Now, arguing  as in Butucea and Taupin \citeyear{buttaupin} we get
that $\vert \mathbb{E}[S_{n} (\theta)] -
S_{\theta^0,P_X}(\theta)\vert^2=o(1)$.

\proof[Proof of the second part of $(\ref{cons})$]
Using that the $Z_i$'s are strictly stationary we get that
\begin{eqnarray*}
\mbox{Var}[S_{n}(\theta)]&=&\mbox{ Var} \left [
n^{-1}\sum_{k=1}^n
\mathbb{R}e\Big[((Z_k-f_\theta)^2\,w)
\star K_{n,C_n}(Z_{k-1})\Big]\right]\\
&\leq &
\frac{1}{n}\mbox{ Var }(A_{1,0})+
\frac{2}{n}\sum_{i=2}^{n}
|\mbox{ Cov}\big( A_{1,0},A_{i,i-1}\big)|\\
&\leq &
\frac{3}{n}\mbox{ Var }(A_{1,0})+
\frac{2}{n}\sum_{k=3}^{n}
|\mbox{ Cov}\big( A_{1,0},A_{k,k-1}\big)|
\end{eqnarray*}
with
$$A_{k,k-1}=\mathbb{R}e\Big[\big((Z_k-f_\theta)^2\,w)
\star K_{n,C_n}(Z_{k-1}) \Big].$$
Arguing  as in Butucea and Taupin \citeyear{buttaupin} we obtain that $\lim_{n \rightarrow \infty}
n^{-1}\mbox{ Var }(A_{1,0})=0$.
It remains to study $$
\frac{1}{n}\sum_{k=3}^{n}
|\mbox{ Cov}\big( A_{1,0},A_{k,k-1}\big)|
 .$$

\begin{lem}
\label{lemcov}
Let $\Psi$ such that $\mathbb{E}(|\Psi(Z)|)<\infty$ and let $\Phi$  be an integrable function. Let
$$B_{k,k-1}=\mathbb{R}\left[e\Psi(Z_k)\Phi\star K_{n,C_n}(Z_{k-1})\right].$$
Then  for $k\geq 3$
\begin{multline*}
\mathrm{ Cov}(B_{k,k-1},B_{1,0})=\mathrm{Cov}[\Psi(Z_k)\Phi\star K_{C_n}(X_{k-1}), \Psi(Z_1)\Phi\star K_{C_n}(X_{0})]\\
=\frac{1}{(2\pi)^2}\iint \Phi^*(t)\Phi^*(s)\mathrm{Cov}\big(\Psi(Z_k)e^{-itX_{k-1}},\Psi(Z_1)e^{-isX_0}\big)
K_{C_n}^*(t)K_{C_n}^*(s)dtds.
\end{multline*}
\end{lem}
\proof[Proof of Lemma \ref{lemcov}:]
By stationarity we write
\begin{eqnarray*}
\mbox{Cov}\big(B_{k,k-1},B_{1,0}\big)=\mathbb{E}(B_{k,k-1}B_{1,0})-\mathbb{E}(B_{k,k-1})\mathbb{E}(B_{1,0})=\mathbb{E}(B_{k,k-1}B_{1,0})-\big(\mathbb{E}(B_{1,0})\big)^2.
\end{eqnarray*}
Now, we use that the sequences $(X_k)_{k\in \mathbb{Z}}$ and $(\varepsilon_k)_{k\in \mathbb{Z}}$ are independent. This implies that $(Z_1,X_{0})$ is independent of $\varepsilon_{0}$ and thus
\begin{eqnarray*}
\mathbb{E}(B_{1,0})=\frac{1}{2\pi}\mathbb{R}e\int \Phi^*(t)\mathbb{E}[\Psi(Z_1)e^{-itZ_{0}}]\frac{K_{C_n}^*(t)}{f_\varepsilon^*(-t)}dt=\frac{1}{2\pi}\int \Phi^*(t)\mathbb{E}[\Psi(Z_1)e^{-itX_{0}}]K_{C_n}^*(t)dt.
\end{eqnarray*}
In the same way, for $k \geq 3$,
\begin{multline*}
\mathbb{E}(B_{k,k-1}B_{1,0})\\=\frac{1}{(2\pi)^2} \mathbb{E} \iint \Phi^*(s)\Phi^*(t)\Psi(Z_k)\Psi(Z_1)\mathbb{R}e\Big(e^{-itZ_{k-1}}\frac{K_{C_n}^*(t)}{f_\varepsilon^*(-t)}\Big)
\mathbb{R}e\Big(e^{-isZ_0}\frac{K_{C_n}^*(s)}{f_\varepsilon^*(-s)}\Big)
dtds\\
=\frac{1}{(2\pi)^2}\iint \Phi^*(s)\Phi^*(t)\mathbb{E}(\Psi(Z_k)e^{-itX_{k-1}}\Psi(Z_1)e^{-isX_0})
K_{C_n}^*(t)K_{C_n}^*(s)dtds\, ,
\end{multline*}
and the lemma is proved. \hfill $\Box$

\medskip

It follows from Lemma \ref{lemcov} that for $k\geq 3$, $$\mbox{ Cov}\big( A_{k,k-1},A_{1,0}\big)=
\mbox{Cov}\Big[ \big((Z_k-f_\theta)^2\,w\big)
\star K_{C_n}(X_{k-1}), \big((Z_1-f_\theta)^2\,w\big)
\star K_{C_n}(X_{0})
 \Big]= \sum_{i=1}^9 C_i,
$$
with
\begin{eqnarray*}
C_1&=&\frac{1}{(2\pi)^2}\iint \mbox{Cov}(e^{-itX_{k-1}},e^{-isX_0}) (w f_\theta^2)^*(t)(w f_\theta^2)^*(s) K_{C_n}^*(s)K_{C_n}^*(t)dtds,\\
C_2&=&\frac{1}{\pi^2}\iint \mbox{Cov}(X_ke^{-itX_{k-1}},X_1e^{-isX_0}) (w f_\theta)^*(t)(w f_\theta)^*(s) K_{C_n}^*(t)K_{C_n}^*(s)dtds,\\
C_3&=&\frac{1}{(2\pi)^2}\iint \mbox{Cov}[(X_k^2+\varepsilon_k^2)e^{-itX_{k-1}},(X_1^2+\varepsilon_1^2)e^{-isX_0}] w^*(t)w^*(s) K_{C_n}^*(t)K_{C_n}^*(s)dtds,\\
C_4&=&\frac{-1}{2\pi^2} \iint \mbox{Cov}(X_ke^{-itX_{k-1}},e^{-isX_0}) (w f_\theta)^*(t)(w f_\theta^2)^*(s) K_{C_n}^*(s)K_{C_n}^*(t)dtds,\\
C_5&=&\frac{-1}{2\pi^2} \iint \mbox{Cov}(e^{-itX_{k-1}},X_1e^{-isX_0}) (w f_\theta)^*(s)(w f_\theta^2)^*(t) K_{C_n}^*(s)K_{C_n}^*(t)dtds,\\
C_6&=&\frac{1}{(2\pi)^2} \iint \mbox{Cov}[(X_k^2+\varepsilon_k^2)e^{-itX_{k-1}},e^{-isX_0}] w^*(t)(w f_\theta^2)^*(s) K_{C_n}^*(s)K_{C_n}^*(t)dtds,\\
C_7&=&\frac{1}{(2\pi)^2} \iint \mbox{Cov}[e^{-itX_{k-1}},(X_1^2+\varepsilon_1^2)e^{-isX_0}] w^*(s)(w f_\theta^2)^*(t) K_{C_n}^*(s)K_{C_n}^*(t)dtds,\\
C_8&=&\frac{-1}{2\pi^2} \iint \mbox{Cov}[(X_k^2+\varepsilon_k^2)e^{-itX_{k-1}},X_1e^{-isX_0}] w^*(t)(w f_\theta)^*(s) K_{C_n}^*(s)K_{C_n}^*(t)dtds,\\
C_9&=&\frac{-1}{2\pi^2} \iint \mbox{Cov}[X_ke^{-itX_{k-1}},(X_1^2+\varepsilon_1^2)e^{-isX_0}] w^*(s)(w f_\theta)^*(t) K_{C_n}^*(s)K_{C_n}^*(t)dtds
\end{eqnarray*}
Easy computations give
\begin{multline*}
\mbox{Cov}[(X_k^2+\varepsilon_k^2)e^{-itX_{k-1}},(X_1^2+\varepsilon_1^2)e^{-isX_0}]=\\
\mbox{Cov}(X_k^2e^{-itX_{k-1}},X_1^2e^{-isX_0})+\sigma_{\varepsilon}^2\mbox{Cov}(X_k^2e^{-itX_{k-1}},e^{-isX_0})\\+
\sigma_{\varepsilon}^2
\mbox{Cov}(e^{-itX_{k-1}},X_1^2e^{-isX_0})+
\sigma_{\varepsilon}^4\mbox{Cov}(e^{-itX_{k-1}},e^{-isX_0})\,
,\end{multline*}
$$
\mbox{Cov}[(X_k^2+\varepsilon_k^2)e^{-itX_{k-1}},e^{-isX_0}]=\mbox{Cov}(X_k^2e^{-itX_{k-1}},e^{-isX_0})+
\sigma_{\varepsilon}^2\mbox{Cov}(e^{-itX_{k-1}},e^{-isX_0})\, ,
$$
$$
\mbox{Cov}[(X_k^2+\varepsilon_k^2)e^{-itX_{k-1}},X_1e^{-isX_0}]=\mbox{Cov}(X_k^2e^{-itX_{k-1}},X_1e^{-isX_0})
+\sigma_{\varepsilon}^2\mbox{Cov}(e^{-itX_{k-1}},X_1e^{-isX_0})\, .
$$
which induces the decomposition
$\mbox{ Cov}\big( A_{k,k-1},A_{1,0}\big)=  \sum_{i=1}^9 E_i$,
with
\begin{multline*}
E_1=\frac{1}{(2\pi)^2}\iint \mbox{Cov}(e^{-itX_{k-1}},e^{-isX_0})K_{C_n}^*(t)K_{C_n}^*(s) \\
\times [(w f_\theta^2)^*(t)(w f_\theta^2)^*(s)+\sigma_{\varepsilon}^4w^*(t)w^*(s)
+\sigma_{\varepsilon}^2w^*(t)(wf_\theta)^*(s) +\sigma_{\varepsilon}^2w^*(s)(wf_\theta)^*(t))]dtds,
\end{multline*}
\begin{eqnarray*}
E_2&=&C_2=\frac{1}{\pi^2}\iint \mbox{Cov}(X_ke^{-itX_{k-1}},X_1e^{isX_0})(w f_\theta)^*(t)(w f_\theta)^*(s)K_{C_n}^*(t)K_{C_n}^*(s)dtds,\\
E_3&=&\frac{1}{(2\pi)^2}\iint \mbox{Cov}(X_k^2e^{-itX_{k-1}},X_1^2e^{-isX_0})w^*(t)w^*(s)K_{C_n}^*(t)K_{C_n}^*(s)dtds,\\
E_4&=&\frac{-1}{2\pi^2} \iint \mbox{Cov}(X_ke^{-itX_{k-1}},e^{-isX_0})K_{C_n}^*(s)K_{C_n}^*(t) (w f_\theta)^*(t)((w f_\theta^2)^*(s)+\sigma^2_\varepsilon w^*(s)) dtds,\\
E_5&=&\frac{-1}{2\pi^2} \iint \mbox{Cov}(e^{-itX_{k-1}},X_1e^{-isX_0})K_{C_n}^*(s)K_{C_n}^*(t) (w f_\theta)^*(s)((w f_\theta^2)^*(t)+\sigma^2_\varepsilon w^*(t)) dtds,\\
E_6&=&\frac{1}{(2\pi)^2}\iint \mbox{Cov}(X_k^2e^{-itX_{k-1}},e^{-isX_0})K_{C_n}^*(t)K_{C_n}^*(s)w^*(t)
(\sigma_{\varepsilon}^2w^*(s)+(w f_\theta^2)^*(s))dtds,\\
E_7&=& \frac{1}{(2\pi)^2}\iint \mbox{Cov}(e^{-itX_{k-1}},X_1^2e^{isX_0})K_{C_n}^*(t)K_{C_n}^*(s)w^*(s)(\sigma_{\varepsilon}^2w^*(t)+(w f_\theta^2)^*(t)) dtds,\\
E_8&=&\frac{-1}{2\pi^2} \iint \mbox{Cov}(X_k^2e^{-itX_{k-1}},X_1e^{-isX_0}) w^*(t)(w f_\theta)^*(s) K_{C_n}^*(s)K_{C_n}^*(t)dtds,\\
E_9&=&\frac{-1}{2\pi^2} \iint \mbox{Cov}(X_ke^{-itX_{k-1}},X_1^2e^{-isX_0}) w^*(s)(w f_\theta)^*(t) K_{C_n}^*(s)K_{C_n}^*(t)dtds\, .
\end{eqnarray*}
Using
\eref{covalpha} and \eref{Mchains}, we have the upper bounds
\begin{eqnarray*}
|\mbox{Cov}(e^{-itX_{k-1}},e^{-isX_0})|&\leq& C \alpha_{{\bf X}}(k-1) \\
|\mbox{Cov}(X_ke^{-itX_{k-1}},X_1e^{-isX_0})|&\leq& C\int_0^{\alpha_{{\bf X}}(k-2)}Q^2_{|X|}(u) du \\
|\mbox{Cov}(X_k^2e^{-itX_{k-1}},X_1^2e^{-isX_0})|&\leq& C\int_0^{\alpha_{{\bf X}}(k-2)}Q^4_{|X|}(t) dt \\
|\mbox{Cov}(X_k^2e^{-itX_{k-1}},e^{-isX_0})| & \leq & C \int_0^{\alpha_{{\bf X}}(k-1)}Q^2_{|X|}(t) dt\\
|\mbox{Cov}(e^{-itX_{k-1}},X_1^2e^{-isX_0})|     &     \leq     &     C
\int_0^{\alpha_{{\bf X}}(k-2)}Q^2_{|X|}(t) dt\\
|\mbox{Cov}(X_k^2e^{-itX_{k-1}},X_1e^{-isX_0})| & \leq & C \int_0^{\alpha_{{\bf X}}(k-2)}Q^3_{|X|}(t) dt\\
|\mbox{Cov}(X_ke^{-itX_{k-1}},X_1^2e^{-isX_0})|     &     \leq     &     C
\int_0^{\alpha_{{\bf X}}(k-2)}Q^3_{|X|}(t) dt\, .
\end{eqnarray*}
Since ${\mathbb E}(X_1^4) < \infty$ and $\lim_{ k \rightarrow \infty}\alpha_{{\bf X}}(k)=0$, we infer that
  $\lim_{k \rightarrow \infty}|\mbox{ Cov}\big( A_{k,k-1},A_{1,0}\big)|=0$. Now, by Cesaro's mean convergence theorem
$$
\lim_{n \rightarrow \infty}\frac{1}{n}\sum_{k=3}^{n}
|\mbox{ Cov}\big( A_{1,0},A_{k,k-1}\big)|=0 \, .
$$
This completes the proof of the consistency.

\subsection{Proof of Theorem \ref{thv1t}}

\proof[Proof of 1) in Theorem \ref{thv1t}]
Starting from the decomposition \eref{base}
we shall check the three following points.
\begin{itemize}
\item[i)] $\mathbb{E}\left[
  (S^{(1)}_{n}(\theta^0)- S^{(1)}_{\theta^0,P_X}(\theta^0))
( S^{(1)}_{n}(\theta^0)-S^{(1)}_{\theta^0,P_X}(\theta^0))^\top\right]=
O[\varphi_n \varphi_n^\top]$
\item[ii)] $S^{(2)}_{n}(\theta^0) \cvp S^{(2)}_{\theta^0,P_X}(\theta^0);$
\item[iii)] $R_{n}$ defined in \eref{Rn} satisfies
  $R_{n}\cvp 0.$
\end{itemize}
The rate of convergence of $\widehat{\theta}$ is thus given by the
order of $$\mathbb{E}\left[
  (S^{(1)}_{n}(\theta^0)-S^{(1)}_{\theta^0,P_X}(\theta^0))
(S^{(1)}_{n}(\theta^0)-S^{(1)}_{\theta^0,P_X}(\theta^0)
)^\top\right].$$

\begin{center}
\textit{\textbf{Proof of i)}}
\end{center}
We first write
\begin{eqnarray*}
 \left(S^{(1)}_{n}(\theta)\right)_i&=&\frac{1}{n}\sum_{k=1}^n \frac{\partial}{\partial \theta_i}
  \mathbb{R}e\left[((Z_k-f_{\theta})^2w)\star K_{n,C_n}(Z_{k-1})-\mathbb{E}[(Z_k-f_{\theta}(X_{k-1}))^2w(X_{k-1})]\right]
\\
&=&\frac{1}{n}\sum_{k=1}^n \left(\frac{\partial}{\partial \theta_i}\mathbb{R}e (Z_k-f_{\theta})^2w\star
K_{n,C_n}(Z_{k-1})-\mathbb{E}\left[\frac{\partial}{\partial \theta_i} (Z_k-f_{\theta}(X_{k-1}))^2w(X_{k-1})\right] \right)\, .
\end{eqnarray*}

\noindent \textbf{Study of the bias.}
As in Butucea and Taupin \citeyear{buttaupin},
we get that
\begin{eqnarray*}
\left|\mathbb{E}\left[\left(S^{(1)}_{n}(\theta^0)\right)_j
\right]\right|\leq C_1(f_{\theta^0},w,
f_\varepsilon)\min\left[B_{n,j}^{[1]}B_{n,j}^{[2]}\right],
\end{eqnarray*}
for $B_{n,j}^{[q]}$, $q=1,2$, defined in
Theorem \ref{thv1t}.

\medskip

\noindent \textbf{Study of the variance.}
For the variance term, note first that
$$
\mbox{Var}\left(\big(
 S^{(1)}_{n}(\theta^0)\big)_j \right) \leq \frac{3}{n} \mbox{Var}(D_{1, 0})
+\frac{2}{n}\sum_{k=3}^n|\mbox{Cov}(D_{1,0},D_{k,k-1})|,
$$
with
$$
\label{Dk}
D_{k,k-1}
=\mathbb{R}e\Big(\big(-2Z_kf^{(1)}_{\theta^0,j}+2f_\theta f^{(1)}_{\theta^0,j}\big)w\Big)\star K_{n,C_n}(Z_{k-1}).
$$
The first part in $\mbox{Var}\big[( S^{(1)}_{n}(\theta^0))_j\big]$
is controlled as in Butucea and Taupin \citeyear{buttaupin} by
\begin{eqnarray}\label{premiercontrol}\frac{1}{n} \mbox{Var}(D_{1, 0})
\leq \frac{C(\sigma_{\xi}^2, f_{\theta^0}, f^{(1)}_{\theta^0,j},w,
f_\varepsilon)}{n}\min\{ V_{n,j}^{[1]}(\theta^0),V_{n,j}^{[2]}(\theta^0)\}
\end{eqnarray}
with $V_{n,j}^{[q]}$, $q=1,2$ defined in  Theorem \ref{thv1t}.
We now control the term
$$\frac{1}{n}\sum_{k=3}^n |\mbox{Cov}(D_{1,0},D_{k,k-1})|.
$$
Applying again Lemma \ref{lemcov}, we obtain that
$$
\mbox{Cov}(D_{1,0},D_{k,k-1}) =F_1+F_2+F_3+F_4
$$
with
\begin{eqnarray*}
F_1&=&\frac{1}{\pi^2}\mathbb{R}e\iint \mbox{Cov}(X_ke^{-itX_{k-1}},X_1e^{-isX_{0}})\big(f^{(1)}_{\theta^0,j}w\big)^*(t)\big(f^{(1)}_{\theta^0,j}w\big)^*(s)K_{C_n}^*(t)K_{C_n}^*(s)dtds\\
F_2&=&\frac{1}{\pi^2}\mathbb{R}e\iint
\mbox{Cov}(e^{-itX_{k-1}},e^{-isX_{0}})\big(f_{\theta^0}
f^{(1)}_{\theta^0,j}w\big)^*(t)\big(f_{\theta^0}
f^{(1)}_{\theta^0,j}
w\big)^*(s)K_{C_n}^*(t)K_{C_n}^*(s)dtds\\
F_3&=&\frac{-1}{\pi^2}\mathbb{R}e\iint \mbox{Cov}(X_ke^{-itX_{k-1}},e^{-isX_{0}})\big(f^{(1)}_{\theta^0,j}w\big)^*(t)\big(f_{\theta^0}f^{(1)}_{\theta^0,j}w\big)^*(s)K_{C_n}^*(t)K_{C_n}^*(s)dtds\\
F_4&=&\frac{-1}{\pi^2}\mathbb{R}e\iint \mbox{Cov}(e^{-itX_{k-1}},X_1e^{-isX_{0}})\big(f^{(1)}_{\theta^0,j}w\big)^*(t)\big(f^{(1)}_{\theta^0,j}w\big)^*(s)K_{C_n}^*(t)K_{C_n}^*(s)dtds.
\end{eqnarray*}
Using  \eref{covalpha} and \eref{Mchains} we have the upper bounds
\begin{eqnarray*}
|\mbox{Cov}(e^{-itX_{k-1}},e^{isX_{0}})| &\leq & C\alpha_{{\bf X}}(k-1)\\
|\mbox{Cov}(X_ke^{-itX_{k-1}},X_1e^{isX_0})| &\leq& C\int_0^{\alpha_{\bf X}(k-2)}Q^2_{|X|}(u) du
\\
|\mbox{Cov}(X_ke^{itX_{k-1}},e^{isX_0})| &\leq & C \int_0^{\alpha_{\bf X}(k-1)}Q_{|X|}(u) du \\
|\mbox{Cov}(e^{itX_{k-1}},X_1e^{isX_0})| & \leq & C \int_0^{\alpha_{\bf X}(k-2)}Q_{|X|}(u) du .
\end{eqnarray*}
Since ${\mathbb E}(X_1^4) < \infty$, we infer that $Q_{|X|}(u) \leq C u^{-1/4}$, and consequently
all the covariance terms are $O(\sqrt{\alpha_{\bf X}(k)})$. Finally, if
$\sum_{k>0} \sqrt{\alpha_{\bf X}(k)} < \infty$, then
$$\frac{1}{n}\sum_{k=3}^n |\mbox{Cov}(D_{1,0},D_{k,k-1})|  \leq \frac Cn \, .
$$
This, together with \eref{premiercontrol}, implies that
$$
\mbox{Var}\left[\big(
 S^{(1)}_{n}(\theta^0)\big)_j \right] \leq \frac Cn \min\{ V_{n,j}^{[1]}(\theta^0),V_{n,j}^{[2]}(\theta^0)\}\, .
$$

\begin{center}
\textit{\textbf{Proof of ii)}}
\end{center}The proof  of \textbf{ii)}  starts from the  expression of  the second
derivative of the estimation criterion
\begin{multline}
\label{d2}
\;\;\left( S^{(2)}_{n}(\theta)\right)_{j,k}
=\frac{1}{2\pi n}\sum_{\ell=1}^n \mathbb{R}e\int \left(-2Z_\ell \frac{\partial^2}{\partial \theta_j \partial \theta_k} (f_{\theta}
    w)
      +\frac{\partial^2}{\partial \theta_j \partial\theta_k} (f_{\theta}^2w)
\right)^*(t)
      \frac{K_{C_n}^*(t)e^{-itZ_{\ell-1}}}{f_\varepsilon^*(-t)}dt.
      \end{multline}
Following the same lines as for the consistency we prove that$$ S^{(2)}_{n}(\theta^0)\cvp S^{(2)}_{\theta^0,P_X}(\theta^0).$$ \qed

\begin{center}
\textit{\textbf{Proof of iii)}}
\end{center}
The proof of \textbf{iii)}  follows from \eref{d2}, from the smoothness
properties of $wf_\theta$ and from Assumption \eref{hypothese}.

\hfill$\Box$

\proof[Proof of 2) in Theorem \ref{thv1t}]
The proof of  2) in theorem \ref{thv1t} is quite  similar to the proof
of 1).  The main  differences appear  in the control  of the  bias and
variance
of $S_n^{(1)}(\theta^0).$ More precisely, we start from
\begin{eqnarray*}
 S^{(1)}_{n}(\theta)
=\frac{1}{n}\sum_{k=1}^n \mathbb{R}e\left(\frac{\partial}{\partial \theta} (Z_k-f_{\theta})^2w\right)\star
K_{n,C_n}(Z_{k-1})-\mathbb{E}\left[\frac{\partial}{\partial \theta} (Z_k-f_{\theta}(X_{k-1}))^2w(X_{k-1})\right] .\end{eqnarray*}

\textbf{Study of the bias}
Since $P_{Z,X}(z,z)=P_X(x)f_\varepsilon(z-x)$ we obtain
that $\mathbb{E}[S_{n}^{(1)}(\theta^0)]-S^{(1)}_{\theta^0,P_X}(\theta^0)$ is equal to
\begin{multline*}
-2 \mathbb{E} \left[f_{\theta^0}(X_0)(f^{(1)}_{\theta^0}
w) \star
K_{C_n}(X_0) - f_{\theta^0} (X_0)f^{(1)}_{\theta^0}(X_0)
w(X_0)\right]
\\+2\mathbb{E} \left[ (f^{(1)}_{\theta^0}
f_{\theta^0}w) \star
K_{C_n}(X_0) -(f^{(1)}_{\theta^0}
f_{\theta^0}w) (X_0)\right],
\end{multline*}
that                                                                 is
$\mathbb{E}[S_{n}^{(1)}(\theta^0)]-S^{(1)}_{\theta^0,P_X}(\theta^0)$
is equal to
\begin{multline*}
-2 \mathbb{R}e\iint
f_{\theta^0}(x)e^{-iux}(f^{(1)}_{\theta^0}w)^*(u)(K_{C_n}^*(u)-1)P_X(dx)\,
du
\\+2 \mathbb{R}e\iint
e^{-iux}(f_{\theta^0}f^{(1)}_{\theta^0}w)^*(u)(K_{C_n}^*(u)-1)P_X(dx)\,
du.
\end{multline*}
It follows that for $j=1,\cdots,d$,
\begin{multline*}
\left\vert \mathbb{E}[(S_{n}^{(1)}(\theta^0))_j] -(S^{(1)}_{\theta^0,P_X}(\theta^0))_j\right\vert \\\leq
\mathbb{E}\vert f_{\theta^0}(X_0)\vert \int \vert (f^{(1)}_{\theta^0,j}w)^*(u)(K_{C_n}^*(u)-1)
\vert du+\int \vert (f_{\theta^0}f^{(1)}_{\theta^0,j}w)^*(u)(K_{C_n}^*(u)-1)
\vert du.
\end{multline*}

\textbf{Study of the variance}
For the study of the variance we combine the proof in Butucea and
Taupin \citeyear{buttaupin} and the proof of \textbf{1)} of Theorem
\ref{thv1t}. For these reasons we only give a sketch of the proof, with
details only for specific parts.
As for the proof of \textbf{1)} we start from
\begin{multline*}
\mbox{Var}\big[(
 S^{(1)}_{n}(\theta^0))_j\big]=\frac{1}{n} \mbox{Var}\left[\mathbb{R}e\left(\frac{\partial [-2Z_k f_\theta
w+f_\theta^2
w ]}{\partial\theta_j}\left.\right|_{\theta=\theta^0}
\right)\star
K_{n,C_n}(Z_{k-1})\right]
\\+\frac{2}{n^2}\sum_{1\leq j<k\leq n}\mbox{Cov}(D_{k,k-1},D_{j,j-1}),
\end{multline*}
with $D_{k,k-1}$ defined in \eref{Dk}.
The control of
$(2/{n^2})\sum_{1\leq j<k\leq n}\mbox{Cov}(D_{k,k-1},D_{j,j-1})$
is done as in the proof of \textbf{1)}.
We now control the first part of $\mbox{Var}\big(( S^{(1)}_{n}(\theta^0))_j\big)$.
\begin{eqnarray*}
\mbox{Var}\big[(S^{(1)}_{n}(\theta^0))_j\big]
\leq\frac{C}{n}\mathbb{R}e\mathbb{E}\left[\left(\frac{\partial [-2Z_i f_\theta
w+f_\theta^2
w ]}{\partial\theta_j}\left.\right|_{\theta=\theta^0}
\right)\star
K_{n,C_n}(Z_i)\right]^2.
\end{eqnarray*}
In other words,
\begin{eqnarray*}
\mbox{Var}\big[(S^{(1)}_{n}(\theta^0))_j\big]
&\leq&\frac{C}{n}\mathbb{R}e\mathbb{E}\left[\left(Z_i f^{(1)}_{\theta^0}
w+f_{\theta^0}f^{(1)}_{\theta^0}w \right)\star
K_{n,C_n}(Z_i)\right]^2\\&=&\frac{C}{n}\mathbb{R}e\mathbb{E}\left[\left((f^2_{\theta^0}(X_0)+\sigma_\xi^2) f^{(1)}_{\theta^0}
w+f_{\theta^0}f^{(1)}_{\theta^0}w \right)\star
K_{n,C_n}(Z_0)\right]^2.
\end{eqnarray*}
Now, write that
\begin{eqnarray*}
\mathbb{R}e\mathbb{E}\left[\left((f^2_{\theta^0}(X_0)+\sigma_\xi^2) f^{(1)}_{\theta^0}
w+f_{\theta^0}f^{(1)}_{\theta^0}w \right)\star
K_{n,C_n}(Z_0)\right]^2=II_1+II_2,
\end{eqnarray*}
with
\begin{eqnarray*}
II_1&=&\mathbb{R}e\iint
f_\varepsilon(z-x)(f_{\theta^0}^2(x)+\sigma_\xi^2)
\left(\int (f^{(1)}_{\theta^0}w)(u) K_{n,C_n}(z-u) du\right)^2 P_{X}(dx) dz
\\II_2&=&\mathbb{R}e\iint
f_\varepsilon(z-x)
\left(\int     (f_{\theta^0}f^{(1)}_{\theta^0}w)(u)     K_{n,C_n}(z-u)
  du\right)^2 P_{X}(dx) dz.
\end{eqnarray*}
We apply H\"{o}lder Inequality and obtain that
\begin{eqnarray*}
\vert II_1 \vert\leq\sup_{z\in
  \mathbb{R}}\mathbb{E}[(f_{\theta^0}^2(X_0)+\sigma_\xi^2)f_\varepsilon(z-X_0) ]
\parallel  (f^{(1)}_{\theta^0}
w)\star
K_{n,C_n}\parallel_2^2,
\end{eqnarray*}
and that $\vert II_1\vert$ is also less than
\begin{eqnarray*}
\mathbb{E}[(f_{\theta^0}^2(X_0)+\sigma_\xi^2) ]
\parallel  (f^{(1)}_{\theta^0}
w)\star
K_{n,C_n}\parallel_\infty^2
\end{eqnarray*}
In the same way we have
\begin{eqnarray*}
\vert II_2\vert\leq \sup_{z\in
  \mathbb{R}}\mathbb{E}[f_\varepsilon(z-X_0) ]
\parallel  (f_{\theta^0}f^{(1)}_{\theta^0}
w)\star
K_{n,C_n}\parallel_2^2,
\mbox{ and }II_2\leq
\parallel  (f_{\theta^0}f^{(1)}_{\theta^0}
w)\star
K_{n,C_n}\parallel_\infty^2.
\end{eqnarray*}
Consequently we have
\begin{eqnarray}
\label{v21}\quad
\mbox{Var}\big[(S^{(1)}_{n}(\theta^0))_j\big]\leq
\frac{C(\sigma_{\xi}^2,f_{\theta^0},f_\varepsilon)}{n}
\left[\parallel  (f^{(1)}_{\theta^0}
w)\star
K_{n,C_n}\parallel_2^2+\parallel  (f_{\theta^0}f^{(1)}_{\theta^0}
w)\star
K_{n,C_n}\parallel_2^2 \right]
\end{eqnarray}
and
\begin{eqnarray}
\label{v22}
\mbox{Var}\big[(S^{(1)}_{n}(\theta^0))_j\big]\leq
\frac{C_1(f_{\theta^0})}{n}
\left[\parallel  (f^{(1)}_{\theta^0}
w)\star
K_{n,C_n}\parallel_2^2+\parallel  (f_{\theta^0}f^{(1)}_{\theta^0}
w)\star
K_{n,C_n}\parallel_1^2 \right].
\end{eqnarray}
By combining \eref{v21} and \eref{v22},
we get that \begin{eqnarray*}
\mbox{Var}\big[(S^{(1)}_{n}(\theta^0))_j\big]\leq \frac{C((f_{\theta^0},\sigma_{\xi}^2,f_\varepsilon)}{n}\min\{ V_{n,j}^{[1]}(\theta^0),V_{n,j}^{[2]}(\theta^0)\}
\end{eqnarray*}
with $V_{n,j}^{[q]}$, $q=1,2$ defined in  Theorem \ref{thv1t}.

\hfill $\Box$

\bibliographystyle{chicago}
\bibliography{Biblio}

\end{document}